\documentclass{nm}
\renewcommand\thisnumber{x}
\renewcommand\thisyear {2025}
\renewcommand\thismonth{January}
\renewcommand\thisvolume{xx}
\usepackage{charter}

\usepackage{multirow}
\usepackage{hyperref}
\hypersetup{
    colorlinks=true,
    linkcolor=black,
    citecolor=black,
    urlcolor=black
}
\usepackage{cleveref} 

\title{A Moving Mesh Isogeometric Method Based on Harmonic Maps}

\author[Tao Wang, Xucheng Meng, Ran Zhang and Guanghui Hu]{Tao Wang\affil{1}, Xucheng Meng\affil{2}\comma\affil{3},
Ran Zhang\affil{1}, Guanghui Hu\affil{4}\comma\affil{5}}
\address{\affilnum{1}\ School of Mathematics, Jilin University, Changchun, China\\
\affilnum{2}\ Research Center for Mathematics, Beijing Normal University, Zhuhai, China\\
\affilnum{3}\ Guangdong Provincial Key Laboratory of Interdisciplinary Research and Application for Data Science, BNU-HKBU United International College, Zhuhai, China\\
\affilnum{4}\ State Key Laboratory of Internet of Things for Smart City and Department of Mathematics, University of Macau, Macao SAR, China\\
\affilnum{5}\ Zhuhai UM Science and Technology Research Institute, Zhuhai, China}
\emails{{\tt wangt21@mails.jlu.edu.cn} (T. Wang), {\tt xcmeng@bnu.edu.cn} (X. Meng), {\tt zhangran@jlu.edu.cn} (R. Zhang), {\tt garyhu@um.edu.mo} (G. Hu)}

\begin{document}

\begin{abstract}
Although the isogeometric analysis has shown its great potential in achieving highly accurate numerical solutions of partial differential equations, its efficiency is the main factor making the method more competitive in practical simulations. In this paper, an integration of isogeometric analysis and a moving mesh method is proposed, providing a competitive approach to resolve the efficiency issue. Focusing on the Poisson equation, the implementation of the algorithm and related numerical analysis are presented in detail, including the numerical discretization of the governing equation utilizing isogeometric analysis, and a mesh redistribution technique developed via harmonic maps. It is found that the isogeometric analysis brings attractive features in the realization of moving mesh method, such as it provides an accurate expression for moving direction of mesh nodes, and allows for more choices for constructing monitor functions. Through a series of numerical experiments, the effectiveness of the proposed method is successfully validated and the potential of the method towards the practical application is also well presented with the simulation of a helium atom in Kohn--Sham density functional theory.
\end{abstract}

\maketitle
\section{Introduction}
In 2005, Hughes et al. \cite{hughes2005isogeometric} proposed a novel numerical method for solving partial differential equations, known as Isogeometric Analysis (IGA), which utilized the same spline basis functions for both geometry representation and solution space construction. The main motivation behind IGA was to exactly represent the geometry and bridge the gap between Computer-Aided Design and Computer-Aided Engineering. Furthermore, in this approach, the method of mesh refinement \cite{hughes2005isogeometric,cottrell2009isogeometric} was simplified as the initial mesh could be constructed based on Non-Uniform Rational B-Splines (NURBS). Similar to the refinements in the traditional finite element method, $h$-refinement, which reduced the mesh size and increased the number of mesh elements, was referred to as \textit{knot insertion} in IGA, and $p$-refinement, which elevated the order of spline basis functions, was known as \textit{order elevation} in IGA. Additionally, a novel refinement approach in IGA, referred to as \textit{$k$-refinement}, introduced in \cite{hughes2005isogeometric} to increase the smoothness of basis functions while significantly reducing the number of degrees of freedom (Dofs).

Compared to the classic finite element methods, the IGA framework can easily provide a high-order and highly regular approximation space, particularly through the $k$-refinement strategy. However, maintaining a balance between accuracy and computational efficiency remains a persistent challenge. Although higher accuracy can be obtained by increasing the number of Dofs, this always results in challenges towards efficiency as a larger linear system of equations is required to solve. To partially resolve this issue, adaptive finite element methods have been developed to enhance the efficiency of numerical methods. Generally speaking, there are three types of adaptive methods, that is, the \textit{$h$-adaptive} method (locally refine the mesh), the \textit{$p$-adaptive} method (locally elevate the degree of basis functions), and the \textit{$r$-adaptive} method (moving mesh method),
and we refer to \cite{Morin_h_refine_2002,WANG20097643,Budd_Huang_Russell_2009,huang2010adaptive} and the references therein for the details.

There have been several attempts in the literature to develop adaptive methods within the framework of IGA. Note that the development of $h$-adaptive techniques for NURBS-based IGA \cite{hughes2005isogeometric} is limited by the tensor-product structure of B-splines in multi-dimensional space \cite{cottrell2009isogeometric}. To overcome this limitation,  several suitable extensions for B-splines have been proposed to admit local mesh refinement. For instance, an adaptive isogeometric method based on hierarchical splines was developed in \cite{buffa2016adaptive}, ensuring the capability of the hierarchical mesh refinement. Additionally, by allowing the hanging nodes, a series of T-splines was developed,  including the analysis-suitable (dual-compatible) T-splines \cite{da2012analysis,scott2012local}, polynomial splines over (hierarchical) T-meshes \cite{deng2006dimensions,deng2008polynomial}, and locally refined splines \cite{johannessen2014isogeometric}. On the other hand, the $p$-adaptive method requires the ability to locally elevate the degree of basis functions. However, the NURBS basis functions \cite{piegl1996nurbs} are constructed with the same degree, limiting the development of \textit{$p$-adaptive} method in IGA, and the hybrid-degree weighted T-splines were introduced by Liu et al. \cite{LIU201642} to make the $p$-adaptive method possible in IGA. Given these considerations, the moving mesh method stands out as a promising approach for improving the computational efficiency of IGA. In fact, the moving mesh method has found widespread applications in various fields, such as computational fluid dynamics \cite{di2005moving,tang2005moving}, phase-field models \cite{wang2008efficient}, and density functional theory \cite{bao2013numerical}. In this work, we aim to develop an integration of IGA and moving mesh method for further enhancing computational efficiency.

It is worthwhile to note the study of $r$-adaptive methods has received
considerable research interests in the IGA community, see, e.g., \cite{BAHARI2024116570,bahari2024,basappa2016,JI2022,ji2023,XU20112021}.
In the previous studies on moving mesh methods, one of the most widely used frameworks, proposed by Winslow \cite{winslow1966numerical}, is to solve a nonlinear, Poisson-like equation for generating a mesh map from logical domain to physical domain \cite{li2001moving}.  There have been several extensions of Winslow's mapping to design the moving mesh methods, see \cite{brackbill1993adaptive,cao1999study}.
Recently, Xu et al. \cite{xu2019efficient} introduced an $r$-adaptive IGA framework based on Winslow's mapping, where the inner control points are updated iteratively. However, the nonlinear Euler–Lagrange equations are needed to solve in the framework of Xu et al. \cite{xu2019efficient},  leading to a high demand on computational cost. 

In this paper, a moving mesh isogeometric method is proposed to numerically solve the PDEs. Within this framework, the numerical solution is obtained by the isogeometric method, and the physical mesh nodes are moved towards the regions where the solution varies rapidly by an iteration method. The moving mesh method used in this paper is based on the harmonic maps, which have been used by Dvinsky \cite{dvinsky1991adaptive} to develop mesh generators. The meshes generated by the  Dvinsky's method are both regular and smooth. Notably, Li et al. \cite{li2001moving,li2002moving} extended the Dvinsky's method and proposed a moving mesh finite element method based on the harmonic mapping, which has been proven to be effective for a wide range of applications, including impressible Naiver--Stokes equations \cite{di2005moving}, diffusion-reaction equations \cite{hu2012moving}, Cahn--Hilliard equation \cite{wang2008efficient} and  Kohn--Sham equation \cite{bao2013numerical}. By integrating this strategy with the IGA framework, we aim to enhance the efficiency of moving mesh method, offering greater applicability for this strategy.

The moving mesh finite element method proposed by Li et al. \cite{li2001moving,li2002moving} requires a weighted average of gradients on a specified set of elements in physical domain, also referred to as the gradient recovery process, to approximate the direction and magnitude of mesh movement of a physical mesh node,  as it uses the piecewise linear $C^0$  finite element method. On the contrary, the spline basis functions can be constructed with $C^k$-regularity ($k\ge1$), and therefore, our moving mesh isogeometric method  can directly and accurately compute the mesh-moving directions and magnitudes of the physical mesh nodes. On the other hand, with the high-order $C^2$ NURBS basis functions, it is possible to accurately and directly use the gradient and high-order derivatives of the numerical solution to construct the monitor function, which is challenging to calculate precisely within the framework of classic $C^0$ finite element method. Moreover, the unique type of $k$-refinement for the isogeometric analysis prevents the explosive increment of the number of Dofs, as demonstrated in \cite{hughes2005isogeometric}. The proposed moving mesh isogeometric method establishes a robust mesh redistribution module that delivers higher accuracy while maintaining a manageable number of Dofs, as the numerical solution algorithm for PDEs and mesh redistribution algorithm are independent of each other, see  \cite{li2001moving,li2002moving} for the details.

The remainder of this paper is structured as follows. \Cref{section2} presents a brief review of isogeometric method. In \Cref{section3}, an algorithm for integrating isogeometric analysis with the moving mesh method is proposed. \Cref{section4} provides a variety of numerical experiments to validate the reliability and effectiveness of the proposed algorithm, along with a practical simulation for a helium atom in the Kohn--Sham density functional theory. Finally, the conclusion is provided in \Cref{section5}.

\section{A brief review of isogeometric analysis}\label{section2}
In this section, a brief review of isogeometric analysis for the Poisson equation will be provided, including the Galerkin formulation, some preliminaries for isogeometric analysis, and \textit{a priori} error estimate derived in \cite{bazilevs2006isogeometric}.

\subsection{The Poisson equation and its isogeometric discretization}

To illustrate the basic idea of isogeometric analysis, we consider the Poisson equation
\begin{equation}\label{PoissonEq}
\left\{
\begin{aligned}
    -\nabla\cdot(\nabla u) &= f \quad \text{in }\Omega,\\
        u & = 0  \quad  \text{on } \partial\Omega,
\end{aligned}\right.
\end{equation}
where $\Omega\subset\mathbb{R}^2$ is a closed bounded domain with Lipschitz continuous boundary $\partial\Omega$, and $f$ is a square integrable function.

The variational formulation of Poisson equation \eqref{PoissonEq} reads:
Find $u\in H_0^1(\Omega)=\{v\in H^1(\Omega):~v|_{\partial\Omega} = 0\}$  such that
\begin{equation}
	\label{PoissonVP}
	a(u,v)=f(v) \quad \forall v\in H_0^1(\Omega),
\end{equation}
where
$$
a(u,v)=\int_{\Omega} \nabla u\cdot \nabla v ~\mathrm{d}x\mathrm{d}y,\quad 
f(v) = \int_{\Omega} f v ~\mathrm{d}x\mathrm{d}y.
$$

The well-posedness of the variational formulation \eqref{PoissonVP} could be achieved by the coercivity and continuity of the bilinear operator $a(\cdot,\cdot)$ over $H_0^1(\Omega) \times H_0^1(\Omega)$, and the Lax-Milgram theorem. In the following, we adopt the IGA method in the framework of Galerkin formulation to numerically discretize \eqref{PoissonVP}.

Given the knot vector
$$\Xi=\{\hat{\xi}_1,\hat{\xi}_2,\cdots,\hat{\xi}_{n+p+1}\},$$ 
where $0=\hat{\xi}_1\leq \hat{\xi}_2\leq \cdots\leq \hat{\xi}_{n+p+1}=1$, $n$ and $p$ are the number and degree of B-spline basis function, respectively, then the $i$-th B-spline basis function of degree $p$, for $1\le i \le n$, denoted by $N_{i,p}(\hat{\xi})$,  is defined by Cox-deBoor recursive formula
\begin{equation}
    \begin{aligned}
&N_{i, 0}(\hat{\xi})= \begin{cases}1, & \text { if } \hat{\xi}_{i} \leq \hat{\xi}<\hat{\xi}_{i+1}, \\
0, & \text { otherwise},\end{cases} \\
&N_{i, k}(\hat{\xi})=\frac{\hat{\xi}-\hat{\xi}_{i}}{\hat{\xi}_{i+k}-\hat{\xi}_{i}} N_{i, k-1}(\hat{\xi})+\frac{\xi_{i+k+1}-\hat{\xi}}{\hat{\xi}_{i+k+1}-\hat{\xi}_{i+1}} N_{i+1, k-1}(\hat{\xi}),
\end{aligned}
\end{equation}
where $k=1,2,\cdots,p$, and we define the quotient $0/0$ as zero. In this work, we use the open knot vector, where the first and last $p+1$ knots are $0$ and $1$, respectively, that is, $\hat{\xi}_1=\cdots=\hat{\xi}_{p+1}=0$ and  $\hat{\xi}_{n+1}=\cdots=\hat{\xi}_{n+p+1}=1$. The B-splines possess attractive properties, such as local support, non-negativity, linearly independent, partition of unity, and global regularity, see, e.g., \cite{hughes2005isogeometric,piegl1996nurbs}. It is noted that the global regularity is associated with the multiplicity of the knot, which describes the number of times a knot appears in the knot vector \cite{piegl1996nurbs}.

The univariate NURBS basis function is defined by 
\begin{equation}
    \label{NURBS1dFormula}
    R_{i,p}(\hat{\xi})=\displaystyle\frac{w_i N_{i,p}(\hat{\xi})}{\displaystyle\sum_{i=1}^n w_i N_{i,p}(\hat{\xi})},\quad i=1,2,\cdots,n,
\end{equation}
where $w_i>0$ is the weight and the B-spline basis function $N_{i,p}(\hat{\xi})$ is defined on a open knot vector $\Xi$.

Given the following two open knot vectors 
\[
\Xi_1=\{0= \hat{\xi}_1,\hat{\xi}_2,\cdots,\hat{\xi}_{n_1+p_1+1}=1\}, 
\text{ and } \Xi_2=\{0=\hat{\eta}_1,\hat{\eta}_2,\cdots,\hat{\eta}_{n_2+p_2+1}=1\},
\]
the partition $\widehat{\mathcal{T}}$ for the parametric domain $\widehat{\Omega} = [0,1]^2$ is defined by
$$\widehat{\mathcal{T}}=\left\{\widehat{K}:~\widehat{K}=[\hat{\xi}_i, \hat{\xi}_{i+1}] \times\left[\hat{\eta}_j, \hat{\eta}_{j+1}\right],~m^*(\widehat{K})\neq 0,~ 1 \leq i \leq n_1, 1 \leq j \leq n_2\right\},$$
where  $m^*(\widehat{K})$ is the measure of the element $\widehat{K}\subset \widehat{\Omega}$. The bivariate NURBS basis functions are defined as
\begin{equation}
    \label{NURBS2dFormula}
    R_{ij}(\hat{\xi},\hat{\eta})=\frac{ w_{i j} N_{i, p}(\hat{\xi}) N_{j, q}(\hat{\eta})  }{\displaystyle\sum_{i=1}^{n_1} \sum_{j=1}^{n_2} N_{i, p}(\hat{\xi}) N_{j, q}(\hat{\eta}) w_{i j}}, \text{ for } 1\le i \le n_1, \,\, 1\le j\le n_2,
\end{equation}
where $N_{i,p}(\hat{\xi})$ and $N_{j,q}(\hat{\eta})$ are the B-spline basis functions in the $\hat{\xi}$- and $\hat{\eta}$-directions, respectively, $p$ and $q$ are the associated degrees of B-splines in $\hat{\xi}$- and $\hat{\eta}$-directions, and $n_1$ and $n_2$ are the corresponding number of basis functions. The bivariate NURBS space $\widehat{\mathcal{N}}_h$ over the parametric domain is given by
$$
\widehat{\mathcal{N}}_h := \operatorname{span}\left\{R_{i j}(\hat{\xi},\hat{\eta}):~1\leq i\leq n_1,1\leq j\leq n_2, ~(\hat{\xi},\hat{\eta})\in \widehat{\Omega}\right\}.
$$
Furthermore, given the control points $\{\mathbf{P}_{i j}\}_{1\le i\le n_1, 1\le j\le n_2}$, the NURBS geometry function $\mathbf{F}$: $\widehat{\Omega}\rightarrow\Omega$ is defined as
\begin{equation}
    \label{GeometryFunc}
    (x,y)=\mathbf{F}(\hat{\xi},\hat{\eta})=\sum_{i=1}^{n_1}\sum_{j=1}^{n_2} R_{i j}(\hat{\xi},\hat{\eta})\mathbf{P}_{i j}~\,\, \forall (\hat{\xi},\hat{\eta})\in \widehat{\Omega}.
\end{equation}
Meanwhile, following the assumption in \cite{bazilevs2006isogeometric}, the geometry function $\mathbf{F}$ is invertible, and has a smooth inverse $\mathbf{F}^{-1}$. With the NURBS geometry function 
$\mathbf{F}:\widehat{\Omega}\rightarrow\Omega$, the partition $\mathcal{T}$ for the physical domain $\Omega$ is given by
$$\mathcal{T}=\left\{K:~K=\mathbf{F}(\widehat{K}),~  \forall \widehat{K}\in \widehat{\mathcal{T}}\right\}.$$
and the bivariate NURBS space $\mathcal{N}_h$ over the physical domain $\Omega$ is defined as
$$
\mathcal{N}_h:= \operatorname{span}\{R_{ij}\circ \mathbf{F}^{-1}:~ 1\leq i\leq n_1,1\leq j\leq n_2\},
$$
then the Galerkin formulation for the variational problem \eqref{PoissonVP} reads: Find $u_h\in V_h$ such that
\begin{equation}
    \label{DiscretePoissonVP}
    a(u_h,v_h)=f(v_h) \quad\forall v_h\in V_h,
\end{equation}
where $V_h= \mathcal{N}_h \cap H_0^1(\Omega)$.

\subsection{\textit{A priori} error estimate}
As pointed out in  \cite{bazilevs2006isogeometric},
the key ingredient for deriving the error analysis in the IGA method is to obtain the interpolation error estimates for NURBS in the physical domain since Cea's lemma still holds in IGA. Compared to the piece-wise polynomials used in the conventional finite element analysis, the support of NURBS basis functions is wider. Following \cite{bazilevs2006isogeometric}, we first define the support extension.

\begin{definition}[Support extension] 
    For a certain non-empty interval $I_j=(\hat{\xi}_j,\hat{\xi}_{j+1})$, we define the support extension of $I_j$ as $$\widetilde{I}_{j}:=\left(\hat{\xi}_{j-p},~\hat{\xi}_{j+p+1}\right).$$
\end{definition}

Note that $\widetilde{I_j}$ can be viewed as the union of the supports of B-spline basis functions whose support intersects $\widetilde{I}_j$. Furthermore, the multi-dimensional support extension for the element $\widehat{K}\subset \widehat{\mathcal{T}}$, denoted by $\widetilde{K}$, can be obtained by using the tensor product structure of NURBS basis functions. 

In what follows, we always assume that the degrees of B-splines in the $\xi$- and $\eta$-directions are the same, that is, $p=q$. Following \cite{bazilevs2006isogeometric}, the approximation error estimates of NURBS can be summarized as follows.

\begin{theorem}\label{NURBSappro}
   Let $ r$ and $s$ be integers satisfying $0\leq r\leq s\leq p+1$, there exists a positive constant $C_{shape}$ such that
    \begin{equation*}
      \sum_{K\in \mathcal{T}}\left|u-\Pi_{\mathcal{N}_h}u \right|^2_{\mathcal{H}^{r}\left(K\right)} \leq C_{shape}\sum_{K\in \mathcal{T}}h_K^{2(s-r)}\sum_{i=0}^s\|\mathbf{\nabla F}\|_{L^\infty(\boldsymbol{F}^{-1}(\widetilde{K}))} ^{2(i-s)}|u|^2_{H^{i}\left(\widetilde{K}\right)}, \quad \forall u\in H^s(\Omega),
    \end{equation*}
    where $\mathcal{H}^m$ is the bent Sobolev space of order $m$, $h_K$ denotes the diameter of the element $K$, $\Pi_{\mathcal{N}_h}: L^2(\Omega)\rightarrow \mathcal{N}_h$ is the projector, and $C_{shape}$ is a constant only depending on the geometry of the physical domain $\Omega$.
\end{theorem}

Using the Cea's lemma and the approximation property of NURBS space in the physical domain $\Omega$, the following theorem can be obtained.
\begin{theorem}\label{h1erroranalysis}
Suppose $u$ is the exact solution of \eqref{PoissonVP} and $u_h$ is the solution of \eqref{DiscretePoissonVP}, let $s$ be an integer  with $s\geq 1$ and assume that $u\in H^s(\Omega)$, then there exists a positive constant $C$ such that
    \begin{equation}
        \|u-u_h\|_{{H}^1(\Omega)}\leq C h^{s-1}\|u\|_{H^s(\Omega)}, \label{H1error}
    \end{equation}
    where $h = \max\{h_K, K\in\mathcal{T}\}$ is the global mesh size.
\end{theorem}

With the Aubin-Nitsche argument and Theorem \ref{h1erroranalysis}, the $L^2$-norm error estimate could be summarized as follows.
\begin{theorem}\label{L2erroranalysis}
     Let $s$ be an integer with $1\le s \le p+1$, and assume that $u$ is the  solution of \eqref{PoissonVP} satisfying $u\in H^s(\Omega)$.
     If $u_h$ is the  solution of discrete variational formulation \eqref{DiscretePoissonVP}, then  there exists a  positive constant $C$ such that
    \begin{equation}
    \|u-u_h\|_{L^2(\Omega)}\leq C h^s\|u\|_{H^s(\Omega)}.  \label{L2error}
    \end{equation}
\end{theorem}

The error estimate in standard isogeometric analysis presented in the Theorem \ref{h1erroranalysis} and Theorem \ref{L2erroranalysis} shows that the optimal convergence order can be achieved when the analytical solution $u$ is sufficiently smooth, as for the classical finite element method. However, the main difference in the estimation is that the full norm of the exact solution is required on the right side of inequalities in \eqref{H1error} and \eqref{L2error}.

\section{A moving mesh isogeometric method} \label{section3}

Following the framework of moving mesh finite element method \cite{li2001moving,li2002moving} based on the harmonic maps, we propose a moving mesh isogeometric method to efficiently solve the PDEs. The implementation details will be presented in the following subsections.

\subsection{Harmonic maps}\label{sec:harmonicmaps}
Following the works of Dvinsky \cite{dvinsky1991adaptive} and Brackbill \cite{brackbill1993adaptive}, suppose there are two compact Riemannian manifolds $\Omega_c$ and $\Omega$ of dimension $n$ with metric tensors $d_{ij}$ and $r_{\alpha\beta}$ in the local coordinates $\boldsymbol{\xi}$ and $\boldsymbol{x}$. We  define the energy $E(\boldsymbol{\xi})$ for a map $\boldsymbol{\xi}:=\boldsymbol{\xi}(\boldsymbol{x})$ as
\begin{equation}
    \label{HarmonicEnergy}
    E(\boldsymbol{\xi})=\frac{1}{2} \int \sqrt{d} d^{i j} r_{\alpha \beta} \frac{\partial \xi^\alpha}{\partial x^i} \frac{\partial \xi^\beta}{\partial x^j} ~\mathrm{d} \boldsymbol{x},~ i,j=1,2,\cdots, n,
\end{equation}
where $d=\operatorname{det}\left(d_{i j}\right),\left(d_{i j}\right)=\left(d^{i j}\right)^{-1}$, and the standard summation convention is assumed. The harmonic maps $\boldsymbol{\xi}$ are the extreme points of the energy functional \eqref{HarmonicEnergy}, which are also the solutions of the Euler-Lagrange equations
\begin{equation*}
\frac{1}{\sqrt{d}} \frac{\partial}{\partial x^i} \sqrt{d} d^{i j} \frac{\partial \xi^k}{\partial x^j}+d^{i j} \Gamma_{\beta \gamma}^k \frac{\partial \xi^\beta}{\partial x^i} \frac{\partial \xi^\gamma}{\partial x^j}=0,
\end{equation*}
where
$$
\Gamma_{\beta \gamma}^k=\frac{1}{2} r^{k \lambda}\left[\frac{\partial r_{\lambda \beta}}{\partial \xi^\gamma}+\frac{\partial r_{\lambda \gamma}}{\partial \xi^\beta}-\frac{\partial r_{\beta \gamma}}{\partial \xi^\lambda}\right]
$$
is the Christoffel symbol of the second kind. With the finite-dimensional Euclidean metric, $\Gamma_{\beta\gamma}^k=0$, we can obtain the simplified energy
\begin{equation}
    \label{HarmonicSimpleEnergy}
    E(\boldsymbol{\xi})= \frac{1}{2} \sum_k\int_\Omega G^{ij}\frac{\partial \xi^k}{\partial x^i}\frac{\partial\xi^k}{\partial x^j}~\mathrm{d}\boldsymbol{x},
\end{equation}
and the corresponding Euler-Lagrange equations
\begin{equation}
    \label{SimpleELequation}
    \frac{\partial}{\partial x^i}\Big(G^{ij}\frac{\partial\xi^k}{\partial x^j}\Big)=0,
\end{equation}
where $G^{ij}=\sqrt{d}d^{ij}$, and its inverse $M=(G^{ij})^{-1}$ is called \textit{monitor function}. Following the notation in \cite{li2001moving}, we refer to $\Omega_c$ and $\Omega$ as logical domain and physical domain, respectively. 

Using the Hamilton-Schoen-Yau theorem \cite{schoen1978univalent,hamilton2006harmonic}, the existence and uniqueness of the harmonic maps $\boldsymbol{\xi}$ are achieved when Riemannian curvature of $\Omega_c$ is non-positive and its boundary $\partial\Omega_c$ is convex. These two conditions are easily satisfied if we let the logical domain  $\Omega_c$ be a cube. Meanwhile, the harmonic map $\boldsymbol{\xi}$ is a continuous, one-to-one mapping with a continuous inverse, which is differentiable and has a non-zero Jacobian.

\subsection{Framework of MMIGM}
In this subsection, we will describe our moving mesh isogeometric method based on harmonic maps in detail. In the two-dimensional case, there are three domains: parametric domain $\widehat{\Omega}$, logical domain $\Omega_c $, and physical domain $\Omega $, which are connected by
$$
\widehat{\Omega} \stackrel{\mathbf{F}}{\longrightarrow}\Omega \stackrel{\boldsymbol{\xi}}{\longrightarrow}\Omega_c ,
$$
where $\mathbf{F}$ is the NURBS geometric mapping defined in \eqref{GeometryFunc}, and $\boldsymbol{\xi}$ is the harmonic map used for mesh redistribution in the moving mesh methods. With the tensor-product structure of NURBS, the procedures presented below can be easily extended to the three-dimensional case.

Following the moving mesh finite element method with the harmonic maps \cite{li2001moving,li2002moving}, the framework of Moving Mesh Isogeometric Methods (MMIGM) can briefly summarized in the Algorithm \ref{alg:MMIGMframework}.

\begin{algorithm}{(\textit{MMIGM framework})}
    \label{alg:MMIGMframework}
    \begin{itemize}
        \item[] \textbf{Step 1:} Generate the initial mesh, denoted by $\boldsymbol{\xi}^0$, for the logical domain $\Omega_c$ by solving the Poisson equation \eqref{InitialwithBoundary};
        \item[] \textbf{Step 2:} Solve the Euler–Lagrange equation \eqref{ELwithBoundary} to obtain a new mesh for $\Omega_c$, denoted by $\boldsymbol{\xi}^*$,  and compute the $L^{\infty}$ error of $\boldsymbol{\xi}^*  - \boldsymbol{\xi}^0$. If the error is less than the prescribed tolerance, the iteration is done; otherwise, move to \textbf{Step 3};
        \item[] \textbf{Step 3:} Compute the direction and magnitude of movement of physical mesh points by using formula \eqref{MMIGMmovement}, then update the physical mesh points with Eq. \eqref{PhysicalMeshUpdate}, and re-parameterize the NURBS  geometry function $\mathbf{F}$;
        \item[] \textbf{Step 4:} Achieve the numerical solution using isogeometric method and the new physical mesh, then update the monitor function in $E(\boldsymbol{\xi})$ and go to \textbf{Step 2}.
    \end{itemize}
\end{algorithm}

In \cite{li2001moving}, the moving mesh method simultaneously redistributes the interior and boundary nodes. Note that the boundary of the physical domain should be preserved with the geometrical constraints, which implies that the vertices of the physical domain will be mapped to the corresponding vertices of the logical domain, respectively \cite{tang2005moving}. For simplicity, we always consider the following mapping set from $\partial \Omega$ to $\partial\Omega_c$:
$$
\begin{aligned}
    \boldsymbol{K}=\{\boldsymbol{\xi}_b\in C^0(\partial\Omega):~\boldsymbol{\xi}_{b}|_{\Lambda_i} \text{ is a piecewise linear mapping}& \\
    \text{without degeneration of Jacobian}&\},
\end{aligned}
$$
where $\boldsymbol{\xi}_b: \partial\Omega\rightarrow\partial\Omega_c$ and $\Lambda_i$ is the set of boundary edges of the physical domain.

With the Dirichlet boundary mapping $\boldsymbol{\xi}_b\in \boldsymbol{K}$, Li et al. \cite{li2001moving} proposed an optimization problem:
\begin{equation}
    \label{HarmonicPwithBoundary}
    \left\{\begin{array}{l}
\min \displaystyle\sum_k \int_{\Omega} G^{i j} \frac{\partial \xi^k}{\partial x^i} \frac{\partial \xi^k}{\partial x^j}~\mathrm{d} \boldsymbol{x}, \\
\text { s.t. }\left.\boldsymbol{\xi}\right|_{\partial \Omega}=\boldsymbol{\xi}_b,
\end{array}\right.
\end{equation}
and the corresponding Euler-Lagrange equation is
\begin{equation}
    \label{ELequationwithboundary}
    \frac{\partial}{\partial x^i}\Big(G^{ij}\frac{\partial\xi^k}{\partial x^j}\Big)=0,
\end{equation}
with Dirichlet boundary condition $\left.\boldsymbol{\xi}\right|_{\partial \Omega}=\boldsymbol{\xi}_b$.

As aforementioned, we will concentrate on the two-dimensional case, where the physical domain $\Omega$ is a polygonal domain, and the corresponding logical domain is convex. Let $\mathcal{T}$ be the quadrilateral mesh for the physical domain, $E_i\in\mathcal{T}$ be the $i$-th element, and $X_j$ be the $j$-th node of the physical mesh $\mathcal{T}$. Similarly, we denote the corresponding notations in the logical domain by $\mathcal{T}_c, E_{i}^c$ and $A_j$, respectively. Furthermore, the set of nodes in physical domain $\Omega$ and logical domain $\Omega_c$ are denoted by $\mathcal{X}=(X_j)$ and $\mathcal{A} = (A_j)$, respectively.

Let $\mathcal{T}^0$ be the initial quadrilateral mesh  in the physical domain $\Omega$, in order to obtain the initial mesh $\mathcal{T}_c$ with nodes $\mathcal{A}=(A_j)$ in the logical domain $\Omega_c$, we solve the following Poisson equation 
\begin{equation}
    \label{InitialwithBoundary}
    \left\{\begin{array}{l}
-\nabla\cdot(\nabla \boldsymbol{\xi}) = 0, \\
\quad \quad ~\left.\boldsymbol{\xi}\right|_{\partial \Omega}=\boldsymbol{\xi}_b,
\end{array}\right.
\end{equation}
where $\boldsymbol{\xi}_b$ is the Dirichlet boundary conditions. We denote the numerical solution for \eqref{InitialwithBoundary} by $\boldsymbol{\xi}^0$. Notice that when the initial mesh $\mathcal{T}_c = \boldsymbol{\xi}^0 $ for the logical domain is available, it will be fixed and used as a reference mesh throughout the moving mesh iterations.

Next, in \textbf{Step 2} of the \Cref{alg:MMIGMframework}, to minimize the energy functional $E(\boldsymbol{\xi})$ and then get the harmonic maps $\boldsymbol{\xi}^*$, we need to solve the equivalent Euler--Lagrange equation
\begin{equation}
    \label{ELwithBoundary}
    \left\{\begin{array}{l}
-\nabla\cdot\left(\displaystyle \frac{1}{M(\boldsymbol{x})}\nabla\boldsymbol{\xi}\right)=0, \\
\qquad\qquad\quad \quad \left.\boldsymbol{\xi}\right|_{\partial \Omega} = \boldsymbol{\xi}_b,
\end{array}\right.
\end{equation}
where $M(\boldsymbol{x})$ is the monitor function.

When the harmonic map $\boldsymbol{\xi}^*$ is computed, the error between $\boldsymbol{\xi}^*$ and $\boldsymbol{\xi}^0$ is subsequently calculated, where $\boldsymbol{\xi}^*$ and $\boldsymbol{\xi}^0$ are solutions to \eqref{ELwithBoundary} and \eqref{InitialwithBoundary}, respectively. If the $L^{\infty}$ norm of error between $\boldsymbol{\xi}^*$ and $\boldsymbol{\xi}^0$ is less than the given tolerance, that is,
$$
    \|\boldsymbol{\xi}^*-\boldsymbol{\xi}^0\|_\infty<\text{tolerance},
$$
then the iteration for the moving mesh in \Cref{alg:MMIGMframework} is stopped, otherwise the movement $\delta \mathcal{X}=(\delta X_j)$ for each inner physical mesh node $X_j$  will be performed to guarantee the movement of mesh nodes in the physical domain, as described  at the \textbf{Step 3}. Compared to the approximate formula used in \cite{li2001moving,li2002moving} to achieve the moving direction and magnitude of a physical mesh node, MMIGM provides a precise formula to calculate the movement of mesh nodes in the physical domain, which is a distinguished feature of our method, as discussed later.  The  formula for determining the direction and magnitude of the movement for physical mesh points will be discussed in the next subsection.

After achieving the movement $\delta X_j$ for each inner node $X_j$ in the physical domain, the corresponding mesh in $\Omega$ can be updated by
\begin{equation}
    \label{PhysicalMeshUpdate}
    \mathcal{X}^* = \mathcal{X}+\tau \delta \mathcal{X},
\end{equation}
where $\tau\in [0,1]$ is a parameter used to avoid mesh wrapping. Based on the updated physical mesh $\mathcal{X}^*$, we could re-parameterize the NURBS geometry function $\mathbf{F}^*$ from $\widehat{\Omega}$ to $\Omega$. Meanwhile, the boundary of the physical domain should be preserved under the assumptions of harmonic maps.

Finally, let's move to the \textbf{Step 4} in \Cref{alg:MMIGMframework}. For the Poisson model problem, which is a time-independent equation, the update for the numerical solution of PDEs could be obtained by the isogeometric method with the new NURBS geometry function $\mathbf{F}^*$. With the updated numerical solution, the new monitor function should be reset to consider the characteristics of numerical solutions, which will be discussed in the following subsection. Then the new harmonic maps will be generated in \eqref{ELwithBoundary} with the updated monitor function.

\subsection{Features of MMIGM}
In this subsection, we will discuss some unique features of MMIGM, including the movement of the mesh nodes $(X_j)$ in the physical domain $\Omega$, the choice of monitor function, and the feasibility of this framework.

For the moving mesh finite element method developed in \cite{li2001moving,li2002moving}, the piecewise linear polynomial is used to construct the finite element space on a triangular mesh, and the moving-direction of $j$-th inner mesh node  $X_j$ in the physical domain $\Omega$ is achieved by the following formula
\begin{equation}
    \label{movementLi}
    \delta X_j=\frac{\displaystyle\sum_{E\in \mathcal{T}_j}| E|~\left(\partial \boldsymbol{x}/\partial \boldsymbol{\xi})\right|_{\text {in}\,  E} \cdot\delta A_j}{\displaystyle\sum_{E\in\mathcal{T}_j}|E|},
\end{equation}
where $\delta A = \boldsymbol{\xi}^0 - \boldsymbol{\xi}^*$ is the difference between the fixed initial logical mesh and the current logical mesh $\boldsymbol{\xi}^*$, $\mathcal{T}_j$ is a set consists of all elements in the physical domain having  $X_j$ as one of its vertices, $|E|$ is the area of the element $E$, and $(\partial\boldsymbol{x}/\partial\boldsymbol{\xi})|_{E}$ is the constant gradient of the piecewise linear map from $\boldsymbol{\xi}^*$ to $\mathcal{X}$ such that $A_i^* \to X_i$ and satisfies the following system of linear equations
$$
\left(\begin{array}{ll}
\displaystyle\frac{\partial x}{\partial \xi} &\displaystyle \frac{\partial x}{\partial \eta} \\
\displaystyle\frac{\partial y}{\partial \xi} & \displaystyle \frac{\partial y}{\partial \eta}
\end{array}\right){\Big|}_{E}
 \left(\begin{array}{ll}
A_{E_{2}^c}^{1}-A_{E_{1}^c}^{1} & A_{E_{3}^c}^{1}-A_{E_{1}^c}^{1} \\
A_{E_{2}^c}^{2}-A_{E_{1}^c}^{2} & A_{E_{3}^c}^{2}-A_{E_{1}^c}^{2}
\end{array}\right)
=  \left(\begin{array}{ll}
X_{E_2}^1-X_{E_1}^1 & X_{E_3}^1-X_{E_1}^1 \\
X_{E_2}^2-X_{E_1}^2 & X_{E_3}^2-X_{E_1}^2
\end{array}\right) ,
$$
where $(X_{E_{i}}^1,X_{E_{i}}^2)$ is the $i$-th ($1\le i \le 3$) vertex of the element $E$ in the physical domain $\Omega$,  and the $i$-th vertex of the corresponding element $E^c$ is denoted by $(A_{E_{i}^c}^1,A_{E_{i}^c}^2)$ in the logical domain $\Omega_c$.

In the previous works \cite{li2001moving,li2002moving}, the piecewise linear polynomials are used to construct the finite element space, which means that the numerical solution is not differentiable at the physical mesh nodes. Therefore, it is inevitable to use a weighted average of gradients on a set of elements to approximate the displacement of a physical mesh node, see the Eq. \eqref{movementLi}. However, for the framework of MMIGM, since NURBS basis functions of degree $p\ge 2$ are used to solve the Euler-Lagrange equations \eqref{ELwithBoundary}, which ensures the finite-dimensional approximation space is of at least global $C^1$ regularity in $\Omega$, and therefore, the gradient $\partial \boldsymbol{x}/\partial\boldsymbol{\xi}$ can be accurately obtained by the inverse coordinate transformation from the physical domain $\Omega$ to the logical domain $\Omega_c$, which is given by
\begin{equation}
    \label{CoordinateTransformation}
    \frac{\partial(x,y)}{\partial(\xi,\eta)}=\left(\begin{array}{cc}
\displaystyle\frac{\partial x}{\partial \xi} & \displaystyle\frac{\partial x}{\partial \eta} \\
\displaystyle\frac{\partial y}{\partial \xi} & \displaystyle\frac{\partial y}{\partial \eta}
\end{array}\right)=
\frac{1}{J}\left(\begin{array}{cc}
~~~\displaystyle\frac{\partial \eta}{\partial y} & \displaystyle-\frac{\partial \xi}{\partial y} \\
\displaystyle-\frac{\partial \eta}{\partial x} & ~~~\displaystyle\frac{\partial \xi}{\partial x}
\end{array}\right)
,~ J = \frac{\partial \xi}{\partial x}\frac{\partial \eta}{\partial y} - \frac{\partial \xi}{\partial y}\frac{\partial \eta}{\partial x},
\end{equation}
where $J=J(x,y)$ is the non-zero Jacobian of the coordinate transformation from  $\Omega_c$ to $\Omega$. Compared with the approximation formula in \eqref{movementLi}, the movement for the $j$-th node of physical mesh, denoted by $\delta X_j $,  is expressed by
\begin{equation}
    \label{MMIGMmovement}
  \delta X_j = \frac{\partial(x,y)}{\partial(\xi,\eta)}\Big|_{A_j}\cdot\delta A_j,
\end{equation}
where $\delta A_j$ is the movement of $j$-th mesh node in the logical domain $\Omega_c$.

The formulas \eqref{CoordinateTransformation} and \eqref{MMIGMmovement} provide an efficient and accurate approach to obtain the movement of the mesh nodes in the physical domain $\Omega$. One possible reason for using the piecewise linear polynomial basis function in the previous moving mesh methods is that it is computationally efficient, and the construction of $C^1$ finite elements is complicated and requires a lot of Dofs in each element. For example, in the 2D case,  one of the typical $C^1$  finite elements is the Argyris triangle \cite{argyris1968tubatriangleelement21} which uses a complete polynomial of degree five and has 21-Dofs in a triangle. However, in the framework of isogeometric analysis, the \textit{$k$-refinement}  is a distinguished method for mesh refinement compared to the conventional finite element method, it not only provides high regularity basis functions but also needs much fewer Dofs than the $p$-refinement strategy.

In addition, only the gradient of numerical solutions is used in the construction of the monitor function as the piecewise linear polynomial basis function is used in \cite{li2001moving}. In the framework of MMIGM, to maintain the mesh quality in the physical domain, the following monitor function with the information of gradient and high-order derivatives could be considered
\begin{equation}
    \label{monitorFunc}
    M(\boldsymbol{x})=\left(\sqrt{\varepsilon+\alpha\left|\nabla u_h^*(\boldsymbol{x})\right|^2+\beta\left|\nabla^2 u_h^*(\boldsymbol{x})\right|^2}\right)I,
\end{equation}
where $u_h^*(\boldsymbol{x})$ is the updated numerical solution of PDEs, $I$ is the identity matrix, and the parameters $\varepsilon,\alpha,\beta$ are properly chosen to mark the regions for the gradient and high-order derivatives information of the numerical solution, respectively. To get the precise value of gradient and high-order derivatives for monitor function, at least three-order basis functions with global $C^2$-regularity are needed to obtain the numerical solutions. Meanwhile, the computational cost in $k$-refinement is much less than the refinements in the traditional finite element method, which makes the numerical experiments possible. Furthermore, based on our numerical experience, the mesh redistribution in the physical domain is sensitive to the parameters in \eqref{monitorFunc}. To prevent mesh wrapping, we prefer the small value of the parameters of \eqref{PhysicalMeshUpdate} and \eqref{monitorFunc}, as discussed in \cite{li2001moving}.

Furthermore, high-order NURBS basis function space provides the capability to calculate the movements of the mesh nodes accurately and get the precise gradient and high-order derivatives in the monitor function \eqref{monitorFunc}. Therefore, with $k$-refinement, the standard NURBS basis functions with global regularity are constructed as the functional basis, which is a more appropriate basis function space with global regularity in moving mesh methods. Additionally, we state that the two-dimensional MMIGM framework can be directly extended to the three-dimensional case which is employed in the practical simulations of the all-electron Kohn--Sham equation for a helium atom in the numerical experiments.

\section{Numerical experiments}\label{section4}
In this section, the performance of the moving mesh isogeometric method is validated through a series of numerical experiments. The first experiment is to provide a comparison of the convergence of numerical errors obtained by the spline-based basis functions with different regularities. Subsequently, we focus on the examinations of the effectiveness and efficiency of \Cref{alg:MMIGMframework} with different  monitor functions, followed by discussions on the Gibbs phenomenon (spurious oscillations) in the experiments. For practical simulations, the method is used to solve the all-electron Kohn--Sham equation for a helium atom. All the experiments are implemented on a personal laptop with the Intel Core i5-10210U CPU (1.60GHz).
\subsection{Case 1: IGA for the Poisson equation}
In this case, we consider the two-dimensional Poisson equation with the non-homogeneous Dirichlet boundary condition. The exact solution is given by
$$
u(x,y) = \sin x\sin y,\quad (x,y)\in\Omega = [-1,1]^2.
$$

To verify the convergence of the isogeometric method \eqref{DiscretePoissonVP},  we take the bi-cubic NURBS basis functions as an example, and  perform the simulation on a series of uniformly refined meshes. In the simulation, by adjusting the multiplicities of inner knots in the knot vectors, we consider the $hp$-refinement, which constructs a global $C^0$-regularity NURBS basis (i.e. the bi-cubic quadrilateral finite element), as well as the $k$-refinement, which generates NURBS basis  functions of $C^2$-regularity.  The corresponding error in the $L^2$-norm and $H^1$-seminorm obtained by the $k$-refinement and the $hp$-refinement are presented in \Cref{tab:k_refinement} and \Cref{tab:hp_refinement}, respectively.

\begin{table}[h]
    \centering
    \centering
    \caption{$k$-refinement}
    \begin{tabular}{ccccc}
        \hline No. of Dofs & $L^2$ error & Order & $H^1$ error & Order \\
        \hline 
        25 & 6.38e-04 & $-$ & 4.89e-03 & $-$ \\
        49 & 3.78e-05 & 4.0782 & 4.93e-04 & 3.3089 \\
        121 & 2.41e-06 & 3.9723 & 6.29e-05 & 2.9702 \\
        361 & 1.57e-07 & 3.9370 & 8.19e-06 & 2.9415 \\
        1,156 & 1.01e-08 & 3.9589 & 1.05e-06 & 2.9620 \\
        4,489 & 6.42e-10 & 3.9773 & 1.33e-07 & 2.9790 \\
        17,161 &4.04e-11 & 3.9881 & 1.68e-08 & 2.9891 \\
        \hline 
    \end{tabular}
    \label{tab:k_refinement}
\end{table}

\begin{table}[h!]
    \centering
    \caption{$hp$-refinement}
    \begin{tabular}{ccccc}
        \hline No. of Dofs & $L^2$ error & Order & $H^1$ error & Order \\
        \hline 
        49 & 1.72e-04 & $-$ &  2.30e-03 & $-$ \\
        169 & 1.15e-05 & 3.9026 & 3.06e-04 & 2.9199 \\
        625 & 7.56e-07 & 3.9280 & 3.84e-05 & 2.9951 \\
        2,401 & 4.85e-08 & 3.9583 & 4.78e-06 & 3.0055 \\
        9,409 & 3.08e-09 & 3.9781 & 5.96e-07 & 3.0048 \\
        37,249 & 1.94e-10 & 3.9888 & 7.43e-08 & 3.0029 \\
        148,225 & 1.22e-11 & 3.9944 & 9.28e-09 & 3.0016 \\
        \hline
    \end{tabular}
    \label{tab:hp_refinement}
\end{table}

It can be observed from the \Cref{tab:k_refinement} and \Cref{tab:hp_refinement} that the numerical convergence rates agree well with the theoretical ones. 
Furthermore, we can  observe that the Dofs of $k$-refinement are significantly fewer than those in $hp$-refinement within the isogeometric framework, while the errors of both refinements remain comparable. Therefore, the $k$-refinement is advantageous over the $hp$-refinement.

\subsection{Case 2: gradient in the monitor function}\label{case_2}
\begin{figure}[h]
    \centering
    \includegraphics[width=0.85\linewidth]{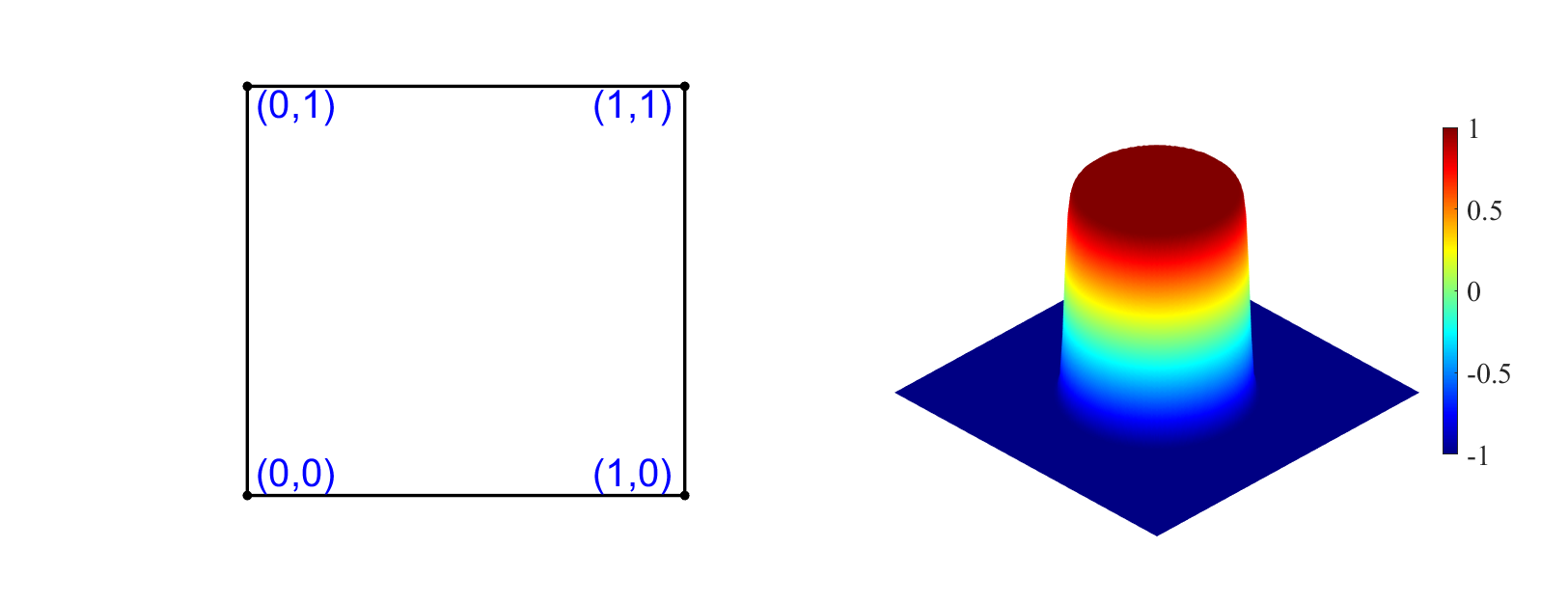}
    \caption{The fixed logical domain (left) and the exact solution \eqref{ExactSolution1} (right).}
    \label{fig:LogicalDomainandExactSolution}
\end{figure}

In this case, we utilize the MMIGM with the gradient of numerical solution as the monitor function to solve the two-dimensional Poisson equation, where the physical domain is $\Omega=[0,1]^2$, and the exact solution is 
\begin{equation}
    \label{ExactSolution1}
    u(x,y) = \tanh{\left(\frac{0.25-\sqrt{(x-0.5)^2+(y-0.5)^2}}{0.01}\right)}.
\end{equation}
The logical domain for MMIGM is a unit square $\Omega_c=[0,1]^2$, as shown in \Cref{fig:LogicalDomainandExactSolution} (left). The exact solution $u$ exhibits large variations along the circle $\mathcal{C}: (x-0.5)^2+(y-0.5)^2= 0.25^2$, as illustrated in  \Cref{fig:LogicalDomainandExactSolution} (right). Therefore, the gradient of the numerical solution $u_h$ obtained  by IGA will be used to construct the monitor function for this example, which is defined as
\begin{equation}\label{gradmonitor}
    M(\boldsymbol{x})=\sqrt{1+\alpha|\nabla u_h(\boldsymbol{x})|^2},
\end{equation}
where $\alpha = 0.1$ is set for avoiding mesh wrapping.

\begin{figure}[h!]
    \centering
    \includegraphics[width=.85\linewidth]{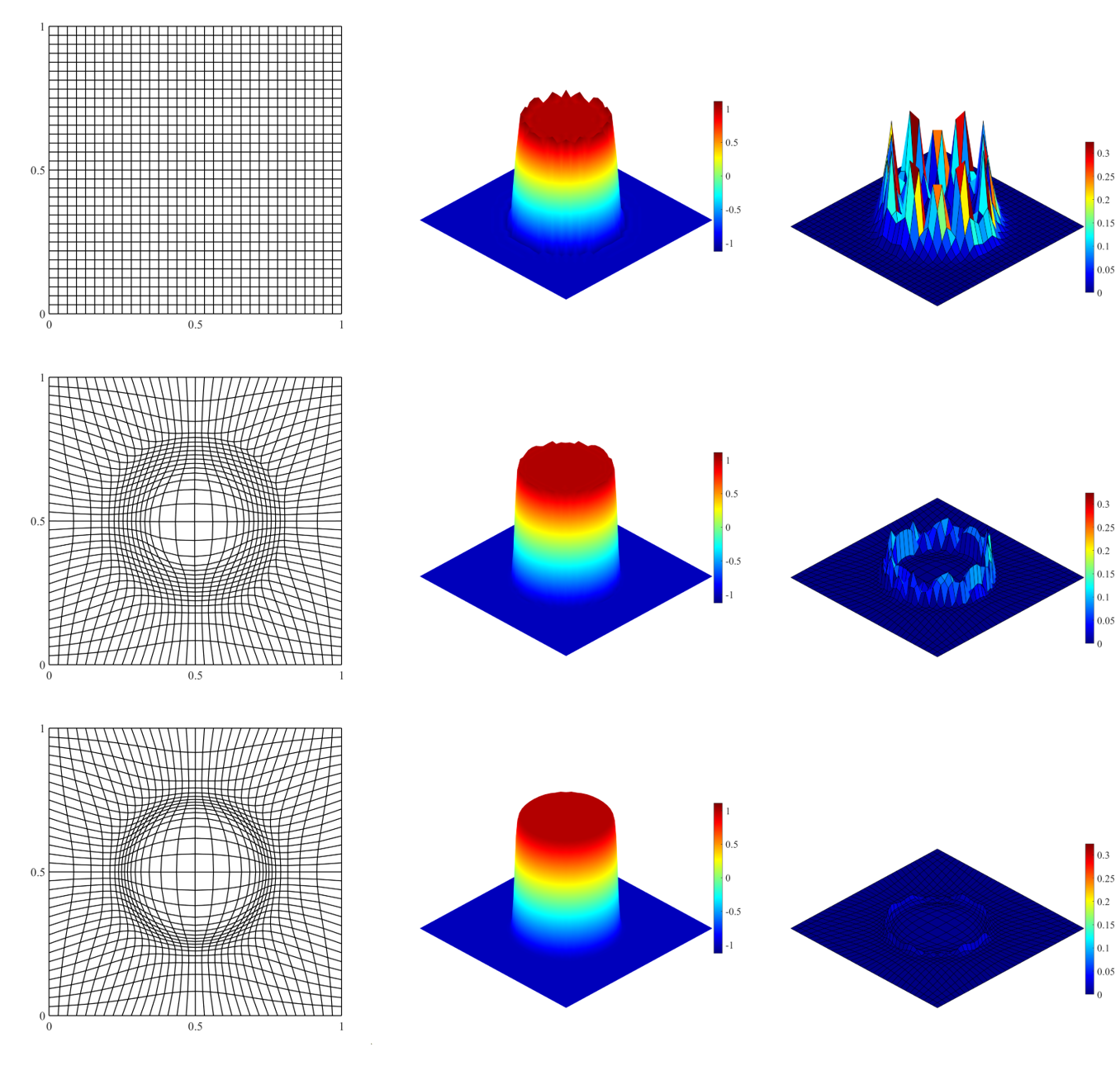}
    \caption{The mesh in the physical domain $\Omega$, corresponding numerical solution and error from left to right at initial (top), intermediate (mid) and final (bottom) iteration with gradients in the monitor function.}
    \label{fig:MMIGMresults}
\end{figure}

Based on the MMIGM framework, the numerical results obtained from the initial iteration, intermediate iteration, and final iteration are presented in \Cref{fig:MMIGMresults}, where we  display the resulting meshes, numerical solutions, and the distribution of $L^2$-error on each element. We observe from  the first column of \Cref{fig:MMIGMresults} that the  physical mesh nodes  are concentrated  at the circle $\mathcal{C}$ where the solution exhibits large variations. Moreover, the numerical oscillations around the circle $\mathcal{C}$ are gradually eliminated through the redistribution of mesh nodes, as shown in the second column of \Cref{fig:MMIGMresults}. Notably, in the final iteration of \Cref{alg:MMIGMframework}, the $L^2$-error in each element decreases to below $0.01$, as shown in the last column of \Cref{fig:MMIGMresults}.

\begin{figure}[h]
    \centering
    \includegraphics[width=0.85\linewidth]{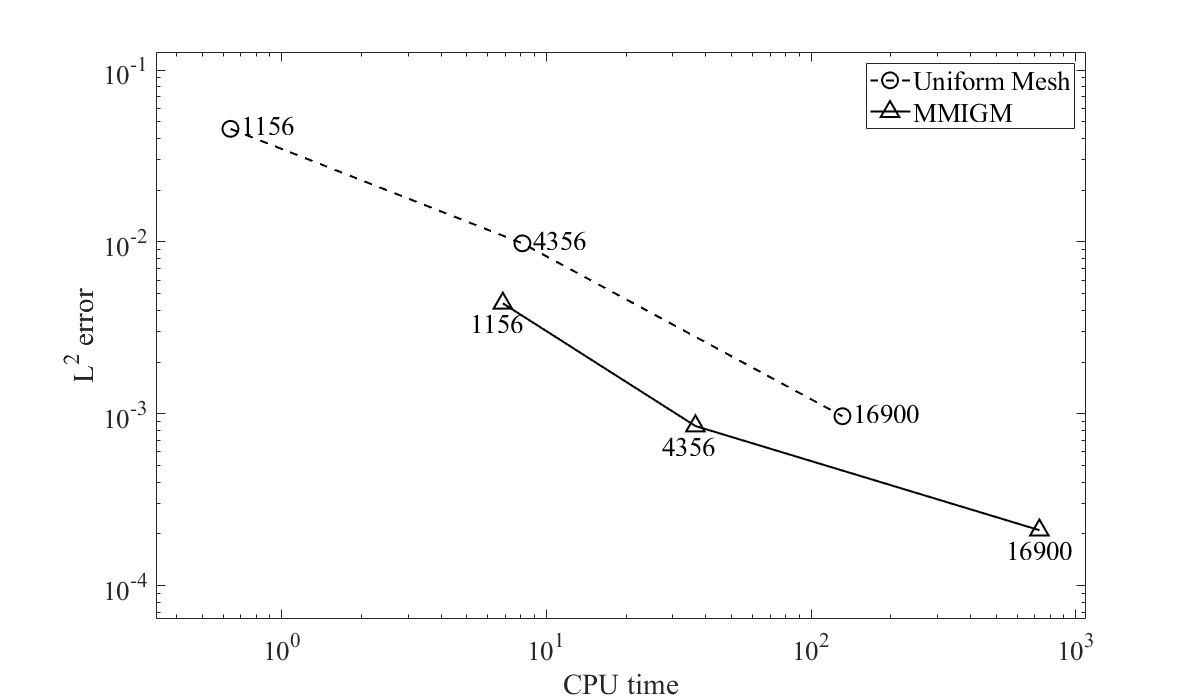}
    \caption{The convergence of the $L^2$ error with respect to the CPU time (s), where the numbers in the figure denote the number of Dofs in each scheme.}
    \label{fig:MMIGM_CPU_error}
\end{figure}

Finally, we compare the numerical results obtained by MMIGM and IGA, and present the related results in \Cref{fig:MMIGM_CPU_error}, where we show the convergence of $L^2$ error with respect to CPU time (in second). It can be observed that the MMIGM can significantly improve the numerical accuracy when compared with the isogeometric method without mesh redistribution. For example, with $16,900$ Dofs, the IGA with uniform mesh can only achieve the error around $10^{-3}$, while after mesh redistribution, the error can be around $10^{-4}$. On the other hand, it can also be found that the MMIGM needs less CPU time to reach the
same error level than the IGA with uniform mesh, which implies that the MMIGM can improve the computational efficiency of IGA.

\subsection{Case 3: high-order derivatives in the monitor function}
In this case, we consider the Poisson equation with the exact solution \eqref{ExactSolution1}, and we utilize the second-order derivatives of numerical solution to construct the monitor function, that is
\begin{equation}
    \label{monitorfunction_sec_derivatives}
    M(\boldsymbol{x})=\sqrt{1+\beta|\nabla^2u_h(\boldsymbol{x})|^2},
\end{equation}
where $\beta = 0.01$ is chosen for avoiding the mesh wrapping.

\begin{figure}[h]
    \centering
    \includegraphics[width=.85\linewidth]{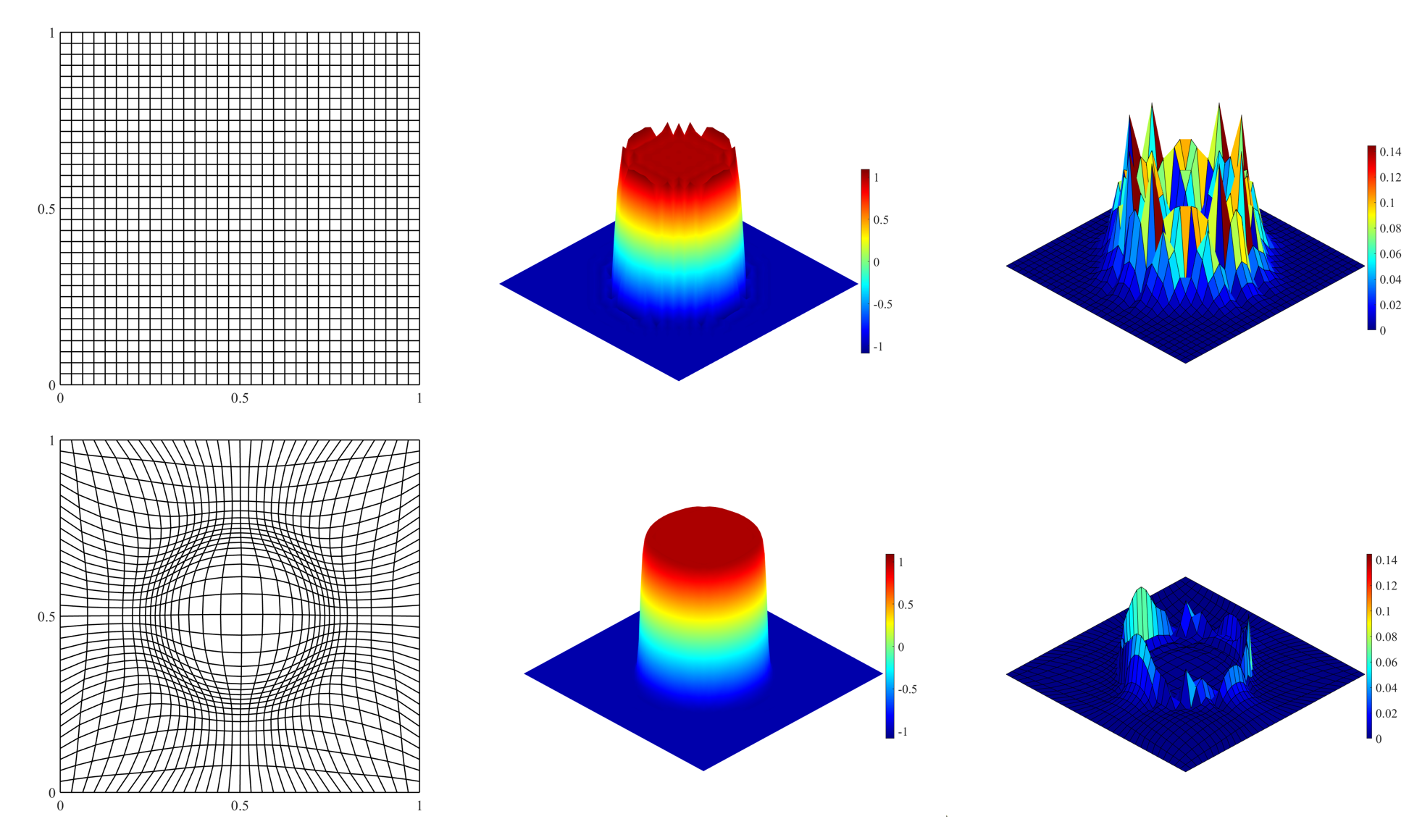}
    \caption{The mesh in the physical domain $\Omega$, corresponding numerical solution and error from left to right at initial (top) and final (bottom) iteration with second-order derivatives in the monitor function.}
    \label{fig:MMIGMresults_deri}
\end{figure}

The initial logical domain for the MMIGM framework is also set as $\Omega_c = [0,1]^2$ and the cubic NURBS basis function space with global $C^2$ regularity is used to directly calculate the second-order derivatives in the monitor function \eqref{monitorfunction_sec_derivatives}. Based on the monitor function introduced in \eqref{monitorfunction_sec_derivatives}, the redistributed mesh, the numerical solutions, and the $L^2$ error for the initial and final iteration are shown in \Cref{fig:MMIGMresults_deri}.
Similar to the numerical results shown in subsection \ref{case_2}, the mesh points are concentrated at the circular curve $\mathcal{C}$ where the solution exhibits large variations, and it be observed that there is no numerical oscillations around $\mathcal{C}$, which successfully demonstrates the effectiveness of MMIGM. Furthermore, as presented in \Cref{fig:MMIGMresults_deri}, the $L^2$ errors in all elements  have been significantly reduced after using the MMIGM, confirming again the effectiveness of our method.

\subsection{Case 4: Gibbs phenomenon at the interface}

For the above two test cases, when the mesh around the circular curve $\mathcal{C}$ is relatively coarse,  we can observe that there exist some numerical oscillations around the curve $\mathcal{C}$, as displayed in  \Cref{fig:MMIGMresults} and \Cref{fig:MMIGMresults_deri}, and this phenomenon is referred to as the Gibbs phenomenon. It is interesting to note that the spurious oscillations can be effectively removed if the MMIGM framework is used.  In this test case, we will provide a detailed discussion on the Gibbs phenomenon, and we consider the Poisson equation with exact solution as in the subsection \ref{case_2}.

\begin{figure}[h]
\centering
\begin{minipage}[h]{0.46\textwidth}     \includegraphics[width=0.85\linewidth]{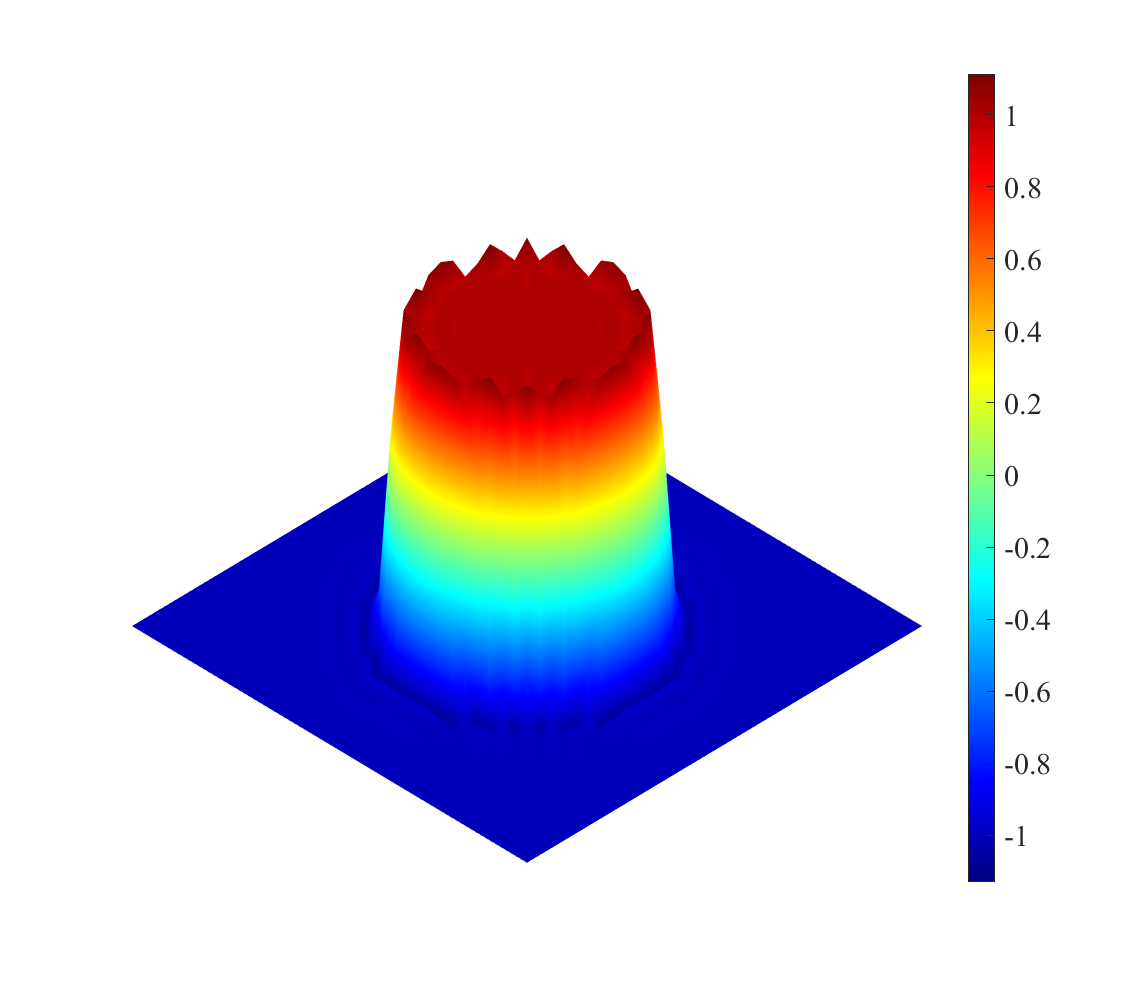}
    \centerline{(a): $k$-refinement ($1156$ Dofs)}
\end{minipage}
\begin{minipage}[h]{0.46\textwidth}    \includegraphics[width=0.85\linewidth]{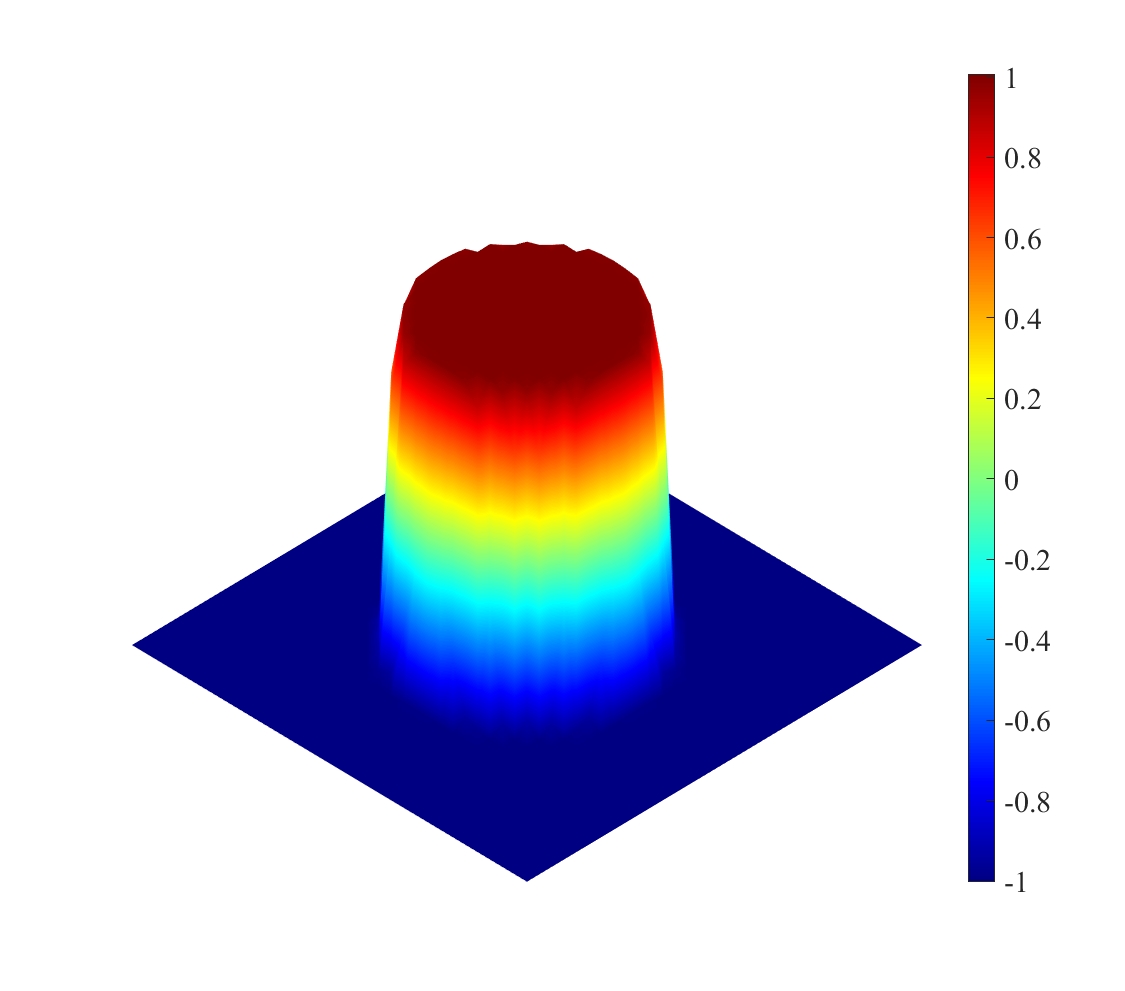}
    \centerline{(b): $hp$-refinement ($4225$ Dofs)}
\end{minipage}
\caption{The numerical solutions obtained by the $k$-refinement, and the $hp$-refinement using IGA with a uniform mesh.}
\label{fig:Gibbs_kandhpcomparsion}
\end{figure}

\begin{figure}[h!]
\centering
\includegraphics[width=0.85\linewidth]{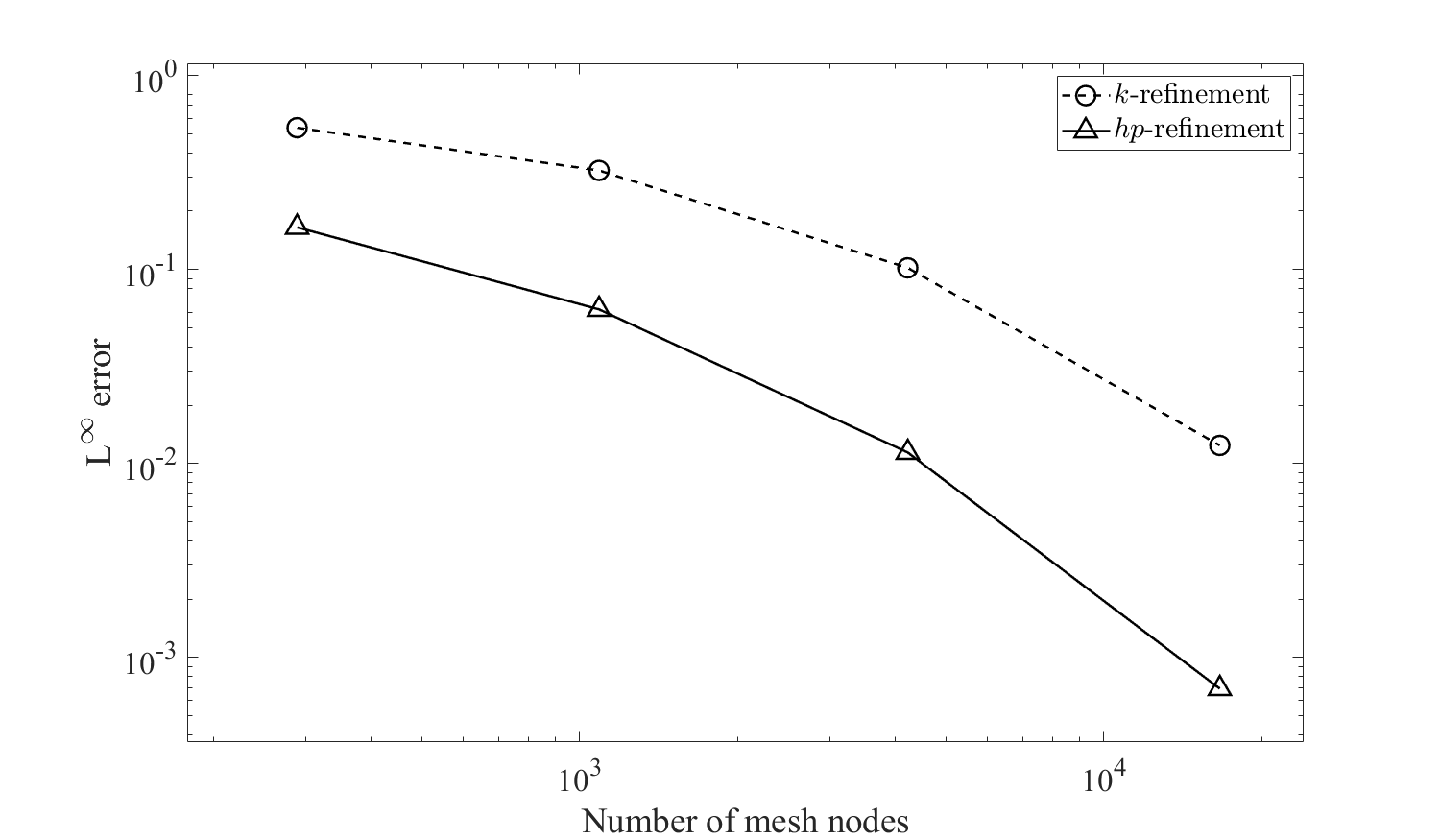}
\caption{The convergence history of $L^\infty$ error associated with the number of mesh nodes, where the IGA with uniform refined meshes is used.}
\label{fig:Gibbserror}
\end{figure}

We employ the isogeometric method and MMIGM with quadratic NURBS basis functions  to compute the numerical solution of this test case. In \Cref{fig:Gibbs_kandhpcomparsion}, we show the numerical results obtained by the isogeometric method using $k$-refinement and $hp$-refinement, where the plotted numerical solutions are obtained on a regular mesh with $32^2$ elements.  It can be observed that 
the spurious oscillations exist for the isogeometric method using both $k$- and $hp$- refinements, and the $k$-refinement will lead to larger oscillations than the $hp$-refinement. One possible explanation is that the $hp$-refinement generates  basis functions with smaller supports, enabling better capability of representing the solution with rapid variations. On the contrary, the $C^{p-1}$ regular NURBS basis functions constructed by $k$-refinement possess wider supports, which may limit their ability to describe the solution with large gradients. The convergence of the $L^{\infty}$ error with respect to the number of mesh nodes is presented in \Cref{fig:Gibbserror}, from which we can see that the smaller error can be achieved by the isogeometric method with the $hp$-refinement.

\begin{figure}[h!]
\centering
\begin{minipage}[t]{0.32\textwidth}
    \includegraphics[width=.9\linewidth]{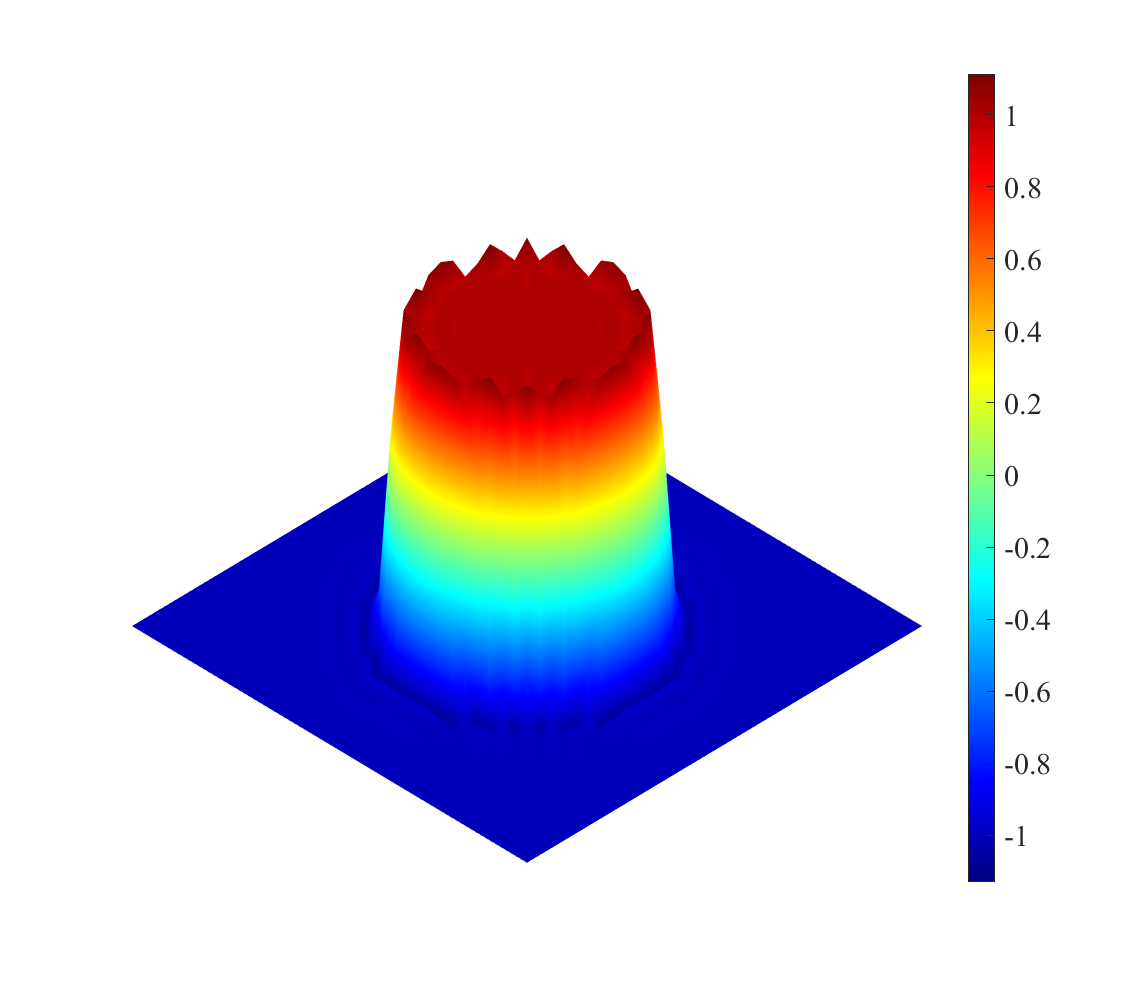}
\end{minipage}
\begin{minipage}[t]{0.32\textwidth}
    \includegraphics[width=.9\linewidth]{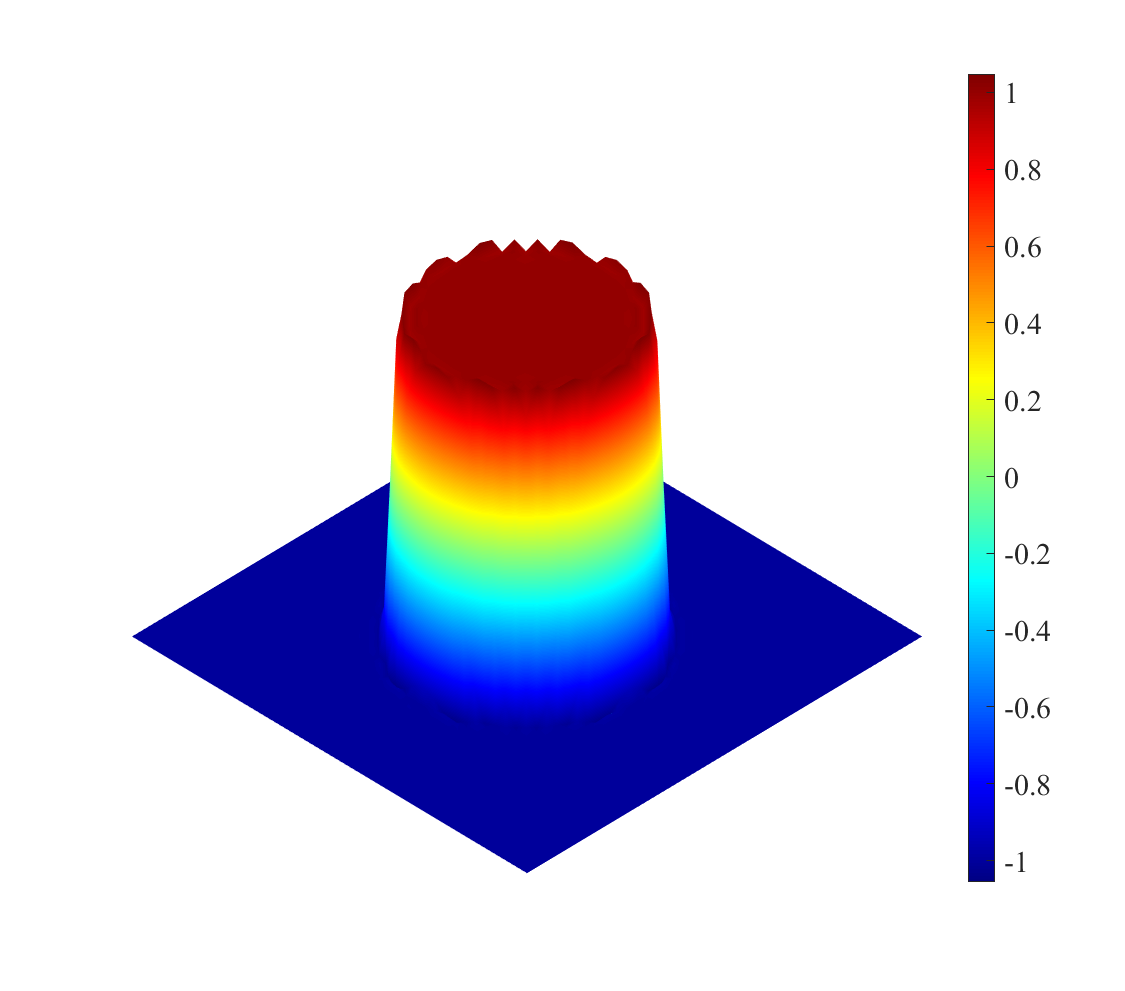}
\end{minipage}
\begin{minipage}[t]{0.32\textwidth}
\includegraphics[width=.9\linewidth]{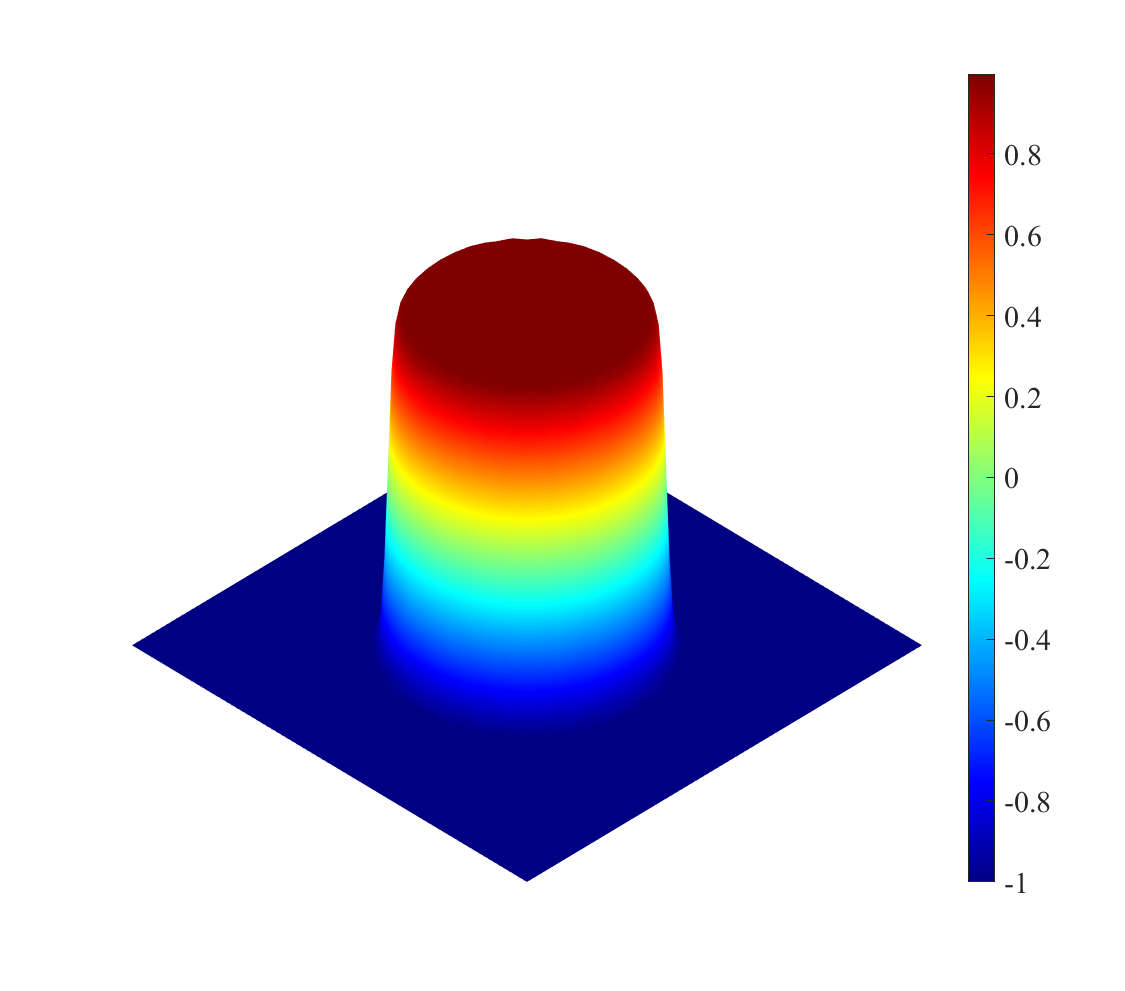}
\end{minipage}
\caption{The numerical solutions of the initial (left), intermediate (middle), and final (right) iteration using $k$-refinement ($1156$ Dofs) and MMIGM framework for the Poisson equation.}
\label{fig:MMIGMforGibbs}
\end{figure}

We next present the numerical results achieved by the MMIGM with $k$-refinement in \Cref{fig:MMIGMforGibbs}, where the  mesh contains $32 \times 32$ elements, and the gradient of numerical solution is used to construct the monitor function. By redistributing the mesh nodes utilizing the MMIGM framework, the numerical oscillations have been effectively removed,  as shown in \Cref{fig:MMIGMforGibbs} (right), which successfully demonstrates the effectiveness of the proposed method to prevent the numerical oscillations.

\subsection{Case 5: the Kohn--Sham equation for a helium atom}
As previously discussed, the MMIGM framework can be straightforwardly extended to the three-dimensional problems. In this test case, the MMIGM framework is employed to simulate the all-electron Kohn--Sham (KS) equation \cite{fiolhais2003primerDFT,bao2012hadaptive} for a helium atom. To handle the singularity arising from the external potential term in the KS equation, the MMIGM framework is utilized to relocate the grids around the nuclear position. There have been several related works \cite{yang2020implicitKSsolver,shen2023allelectron,kuang2024allelectron,li2024method,wang2024HBsplines} devoted to the numerical simulation of all-electron KS equation in recent years. In this context, the self-consistent field (SCF) iterations are used to handle the nonlinearity. The total energy of the ground state is obtained by solving the discretized KS equation generated using NURBS basis functions in the SCF iteration. Furthermore, for the wavefunctions obtained by the SCF iteration, the MMIGM framework is used for the mesh redistribution, providing an optimized initial guess for the subsequent SCF iteration on the $r$-adaptive mesh.

The Kohn--Sham equation could be summarized as
\begin{equation}
    \label{KSequation}
    \left(-\frac12\nabla^2+V_{KS}(\mathbf{r})\right)\psi_i(\mathbf{r}) = \varepsilon_i\psi_i(\mathbf{r}),
\end{equation}
where $-\displaystyle\frac{1}{2}\nabla^2$ denotes the operator of kinetic energy, $(\varepsilon_i,\psi_i(\mathbf{r}))$ denotes the $i$-th eigenpair, and $V_{KS}$ is the effective potential in Kohn--Sham equation and expressed as:
\begin{equation*}
    \label{KSpotential}
    V_{KS}(\mathbf{r}) = -\sum_{j}\frac{Z_j}{|\mathbf{r}-\mathbf{r}_j|} + \int\frac{\rho(\mathbf{r}^\prime)}{|\mathbf{r}-\mathbf{r'}|}~\mathrm{d}\mathbf{r^\prime} +V_{xc}(\mathbf{r},\rho),
\end{equation*}
where $\mathbf{r}_j$ and $Z_j$ are the position and charge of $j$-th nucleus, $\rho(\mathbf{r})=\sum_{i} f_i|\psi_i(\mathbf{r})|^2$ denotes the electron density, and $V_{xc}(\mathbf{r},\rho)$ is the exchange-correlation potential.

\begin{figure}[h]
\centering
    \begin{minipage}[t]{0.45\textwidth}
        \includegraphics[width=.8\linewidth]{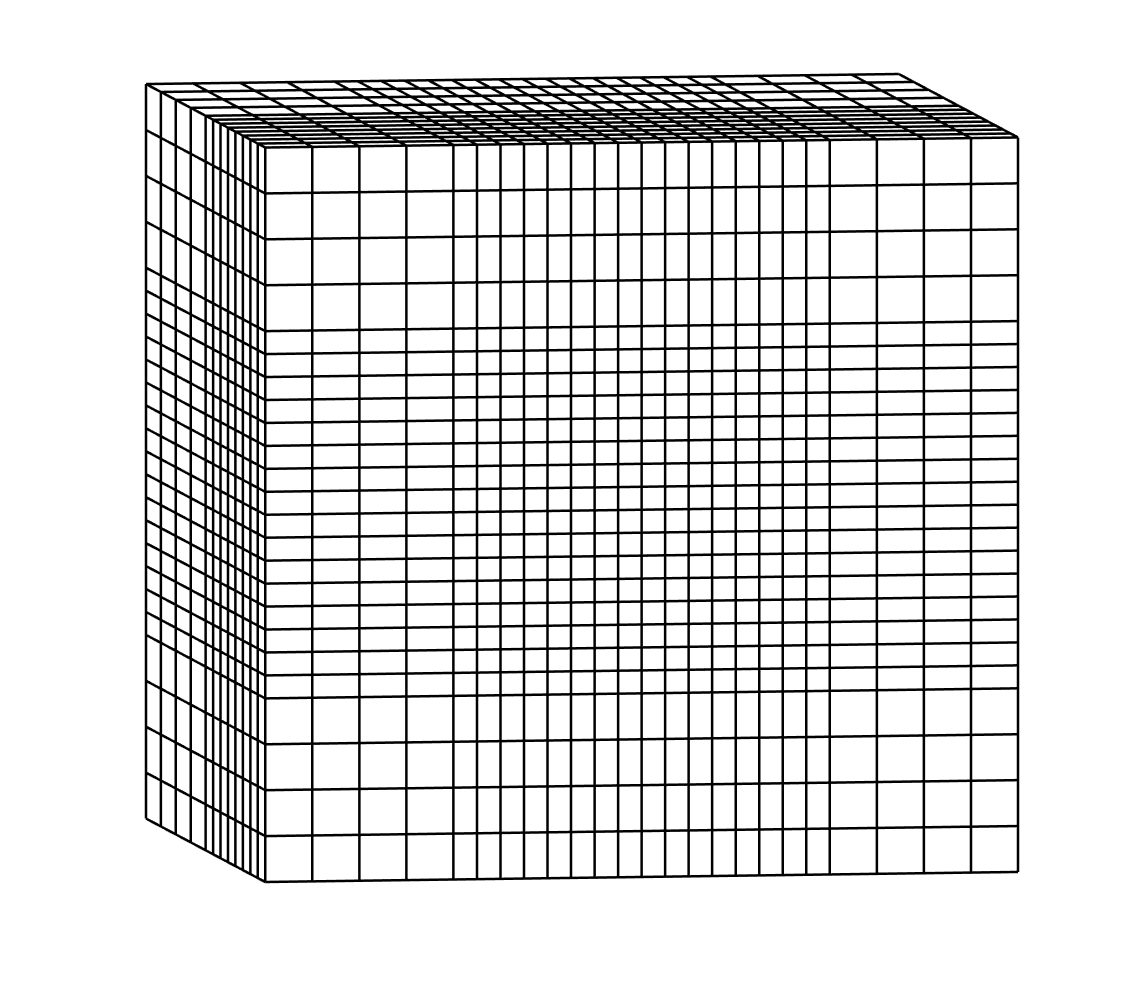}
    \end{minipage}
    \begin{minipage}[t]{0.45\textwidth}
        \includegraphics[width=.8\linewidth]{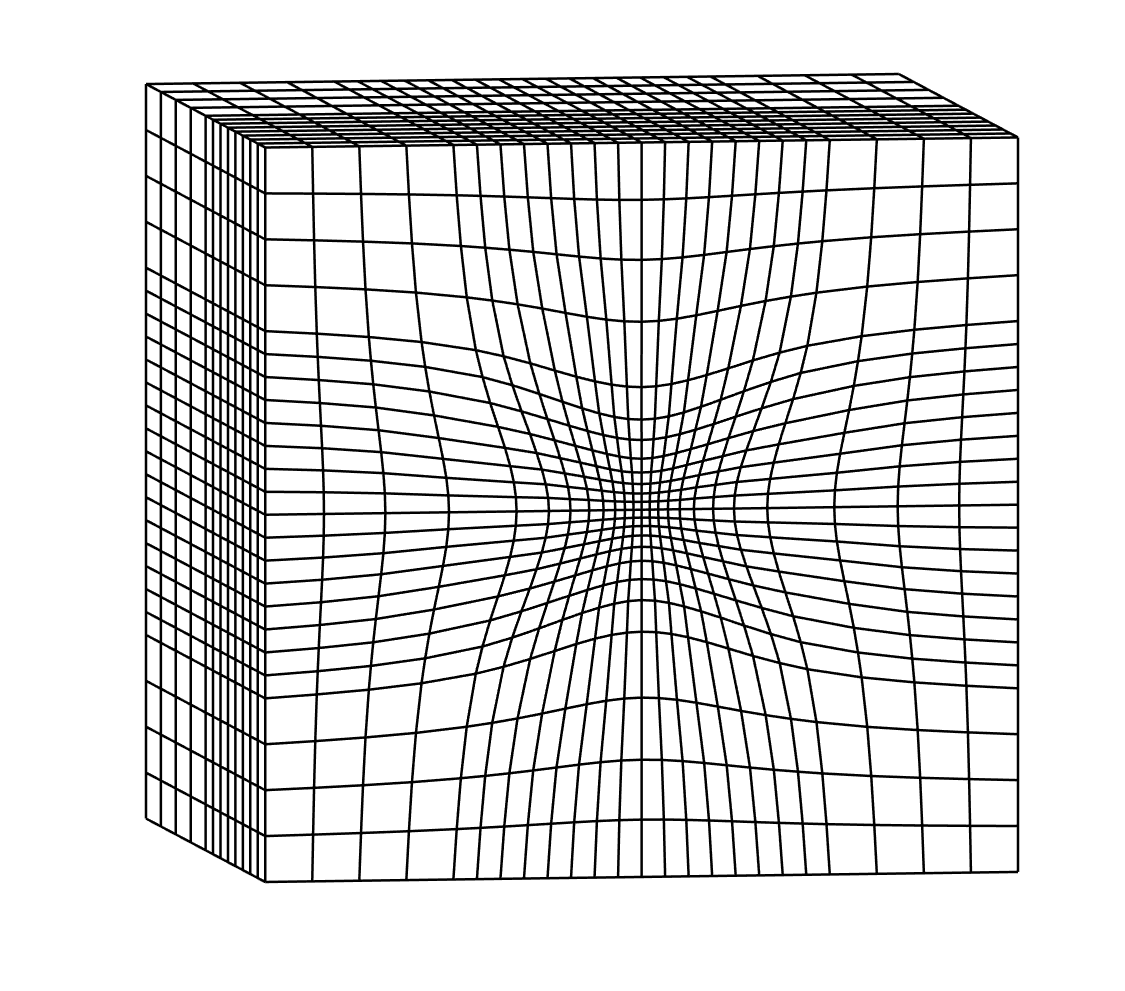}
    \end{minipage}
    \caption{The slice of the three-dimensional mesh of the initial mesh (left) and the redistributed mesh (right) cut with a plane normal to the $z$-axis.}
    \label{fig:MMIGM_Helium_mesh}
    \end{figure}

Here, we will focus on the numerical simulations of the ground state for a helium atom. The Kohn--Sham equation of the helium atom is restricted to the physical domain $\Omega = [-10,10]^3$. Meanwhile, the multipole expansion is employed to generate the boundary conditions of Hartree potential. The reference total energy is $E_{\mathrm{ref}} = -2.8343 \, \mathrm{a.u.}$, which is obtained by the state-of-art package \texttt{NWChem} using the \emph{aug-cc-pv6z} basis set \cite{valiev2010nwchem}. Furthermore, since the singularity of the helium atom is at the nuclei, we prefer a radial mesh as the initial mesh, which distributes more elements around the nucleus. The initial radial mesh is presented in \Cref{fig:MMIGM_Helium_mesh}(a).

\begin{figure}[h]
    \centering
    \includegraphics[width=0.4\linewidth]{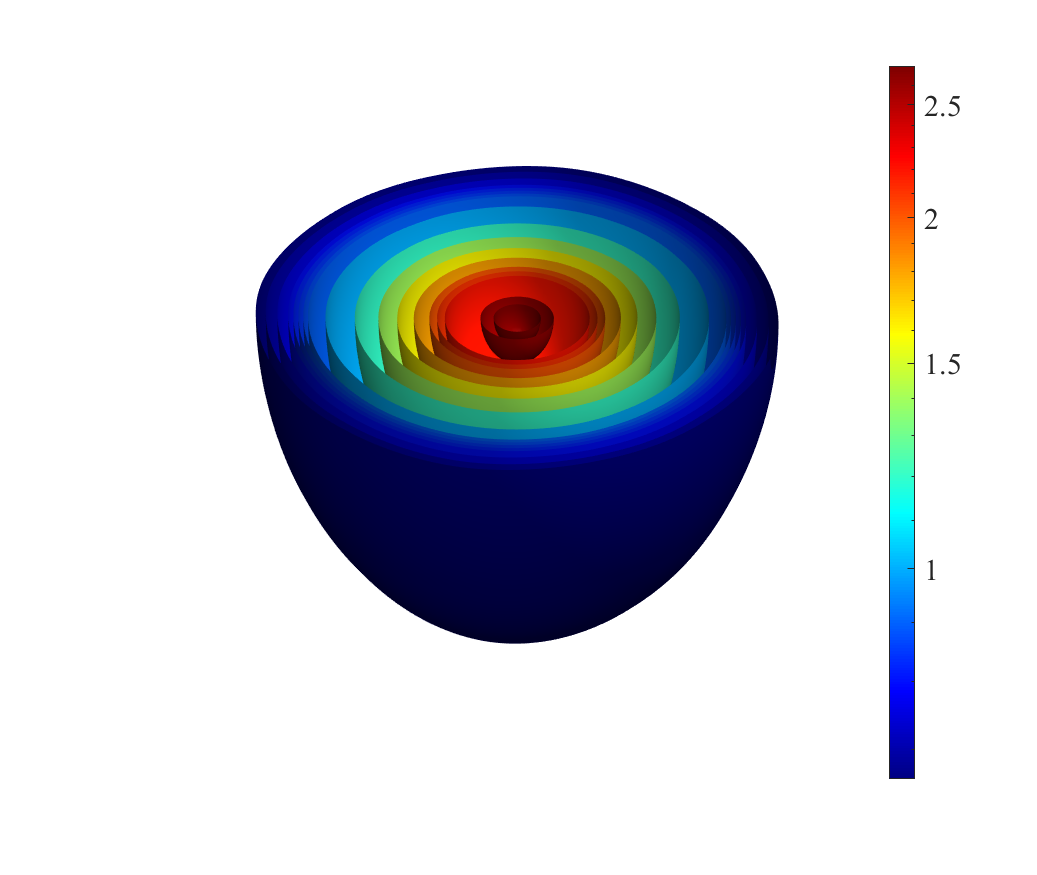}
    \caption{The sliced density of the helium atom simulated by the SCF iteration and MMIGM framework.}
    \label{fig:heliumdensity}
\end{figure}

To obtain the wavefunctions for the KS equation, the locally optimal block preconditioned conjugate gradient (LOBPCG) method \cite{knyazev2007blockLOBPCG} is employed as the eigensolver for the generalized eigenvalue problem arising from the discretized KS equation. Furthermore, a linear mixing scheme for the electron density is adopted with an appropriately chosen mixing parameter $0.618$ for the SCF iterations. The gradient in the monitor function is chosen in the framework of MMIGM. With the combinations of SCF iterations and MMIGM framework, the slice of the three-dimensional mesh shown in \Cref{fig:MMIGM_Helium_mesh}(b) illustrates the redistribution of nodes around the nuclear position. Among the simulations, the ground state energy is improved from $-2.7611~\mathrm{a.u.}$ to $-2.8328~\mathrm{a.u.}$. Based on the $r$-adaptive mesh displayed in \Cref{fig:MMIGM_Helium_mesh}(b), the sliced electron density of the wavefunctions for the helium atom is demonstrated in Figure \ref{fig:heliumdensity}, confirming the capability of our MMIGM framework in the numerical simulations in the three-dimensional space.

\section{Conclusions}\label{section5}
In this paper, the Moving Mesh Isogeometric Method (MMIGM) based on harmonic maps is proposed to improve the computational efficiency of isogeometric method. Thanks to the high regularity of NURBS basis functions, our method can directly and accurately compute the gradient of the map  between current logical mesh and the physical mesh, rather than using an approximation as adopted  in the moving mesh finite element method \cite{li2001moving,li2002moving}. Moreover, the proposed method allows for more choices
for the construction of monitor functions, e.g., the high-order derivatives of numerical solution can be used to construct the monitor function. 

Several numerical experiments are presented to validate the effectiveness of combining the moving mesh method with isogeometric analysis. Firstly, the convergence of IGA based on NURBS basis functions is verified through both $hp$-refinement and $k$-refinement, and the latter requires much fewer Dofs than the former to achieve the same accuracy. The experiments demonstrate that, under the MMIGM framework, the accuracy is improved by incorporating the gradient and high-order derivatives of the numerical solution into the monitor functions. Additionally, we highlight that the numerical oscillations in the Gibbs phenomenon could be eliminated using the MMIGM framework. Besides the two-dimensional Poisson equation, the MMIGM framework is also applied to the practical simulation of the three-dimensional all-electron Kohn--Sham equation, demonstrating an optimized $r$-adaptive mesh for a helium atom. Based on our numerical experience,  we discover that the wavefunctions of the KS equation could be interpreted as an exponential series, given the analytical solution of hydrogen atoms. As a result, we are particularly interested in applying the NURBS basis function for the simulation of the all-electron KS equation, as NURBS basis functions are globally regular over the entire domain. Further details about the numerical simulation for Kohn--Sham equation under the framework of MMIGM will be presented in future work. Additionally, it remains a challenge to apply the MMIGM framework to the time-dependent PDEs \cite{li2001moving}, as a high-order interpolation of $u_h^n$ between two time levels is required. We are committed to exploring this aspect in our future research.

\section*{Acknowledgement}
The authors would like to thank Prof. Tao Tang (Guangzhou Nanfang College \& BNU-HKBU United International College) for his valuable guidance and insightful suggestions on this work.

\bibliographystyle{abbrv}
\bibliography{manuscript}

\end{document}